\documentclass[graybox]{svmult}
\usepackage{amssymb,amsmath,latexsym,amscd,color,stmaryrd}
\usepackage{mathptmx}       
\usepackage{helvet}         
\usepackage{courier}        
\usepackage{makeidx}         
\usepackage{graphicx}        
\usepackage{multicol}        

\makeindex

\def\ga{\mathfrak{a}}

\def\gg{\mathfrak{g}}
\def\gh{\mathfrak{h}}

\def\gk{\mathfrak{k}}
\def\gl{\mathfrak{l}}
\def\gm{\mathfrak{m}}
\def\gn{\mathfrak{n}}

\def\gp{\mathfrak{p}}
\def\gq{\mathfrak{q}}
\def\gr{\mathfrak{r}}
\def\gs{\mathfrak{s}}
\def\gt{\mathfrak{t}}
\def\gu{\mathfrak{u}}
\def\gv{\mathfrak{v}}

\def\gz{\mathfrak{z}}

\def\Ad{{\rm Ad}}
\def\ad{{\rm ad}}
\def\rank{{\rm rank\,}}
\def\trace{{\rm trace\,\,}}
\def\Ind{{\rm Ind\,}}

\def\tr{{\rm trace\,}}
\def\Det{{\rm Det}}
\def\Pf{{\rm Pf}\,}

\def\C{\mathbb{C}}

\def\H{\mathbb{H}}

\def\Q{\mathbb{Q}}
\def\R{\mathbb{R}}

\def\Z{\mathbb{Z}}

\def\cC{\mathcal{C}}
\def\cD{\mathcal{D}}
\def\cE{\mathcal{E}}
\def\cF{\mathcal{F}}

\def\cH{\mathcal{H}}

\def\cL{\mathcal{L}}

\def\cO{\mathcal{O}}

\def\cU{\mathcal{U}}

\title*{Stepwise Square Integrable Representations: \\
the Concept and Some Consequences}

\author{Joseph A. Wolf} 

\institute{Joseph A. Wolf \at Department of Mathematics,
        University of California,
        Berkeley, CA 94720--3840, USA,
	\email{jawolf@math.berkeley.edu}}

\begin{document}
\titlerunning{Stepwise Square Integrable Representations}
\maketitle

\abstract{ There are some new developments on Plancherel formula and 
growth of matrix coefficients for unitary representations of nilpotent 
Lie groups.  These have several consequences for the geometry of weakly 
symmetric spaces and analysis on parabolic subgroups of real semisimple 
Lie groups, and to (infinite dimensional) 
locally nilpotent Lie groups.  Many of these consequences are still
under development.  In this note I'll survey a few of these new aspects
of representation theory for nilpotent Lie groups and parabolic subgroups.}

\section{\hskip -4pt . \hskip 2pt Introduction.} \label{sec1}
\setcounter{lemma}{0}
\setcounter{theorem}{0}
\setcounter{proposition}{0}
\setcounter{corollary}{0}
\setcounter{definition}{0}
\setcounter{remark}{0}
\setcounter{example}{0}

There is a well developed theory of square integrable representations
of nilpotent Lie groups \cite{MW1973}.  It is based on the general 
representation theory of Kirillov \cite{K1962} for connected nilpotent 
real Lie groups.  A connected simply connected Lie group $N$
with center $Z$ is called {\em square integrable} if it has unitary
representations $\pi$ whose coefficients $f_{u,v}(x) =
\langle u, \pi(x)v\rangle$ satisfy $|f_{u,v}| \in \cL^2(N/Z)$.  If $N$ has one
such square integrable representation then there is a certain polynomial
function $P(\gamma)$ on the linear dual space $\gz^*$ of the Lie algebra of
$Z$ that is key to harmonic analysis on $N$.  Here $P(\gamma)$ is the
Pfaffian of the antisymmetric bilinear form on $\gn / \gz$ given by
$b_\lambda(x,y) = \lambda([x,y])$ where $\gamma = \lambda|_{\gz}$\,.  The 
square integrable
representations of $N$ are certain easily--constructed representations
$\pi_\gamma$ where $\gamma \in \gz^*$ with $P(\gamma) \ne 0$,
Plancherel  almost irreducible unitary representations of $N$ are square
integrable, and up to an explicit constant
$|P(\gamma)|$ is the Plancherel density of the unitary
dual $\widehat{N}$ at $\pi_\lambda$.  
This theory has some interesting analytic consequences \cite{W2007}.
\medskip

More recently there was a serious extension of that theory \cite{W2013}.  
Under certain
conditions, the nilpotent Lie group $N$ has a decomposition into
subgroups that have square integrable representations, and the Plancherel
formula then is synthesized explicitly in terms of the Plancherel formulae
of those subgroups.  In particular the extended theory applies to nilradicals
of minimal parabolic subgroups \cite{W2013}.  With a minor technical 
adjustment it has just been extended to nilradicals of arbitrary real
parabolics \cite{W2015c}.  The consequences include explicit Plancherel
and Fourier inversion formulas.  Applications include analysis on minimal
parabolic subgroups \cite{W2014} and, more generally, on maximal amenable
subgroups of parabolics \cite{W2015c},  They also include analysis on
commutative spaces, i.e. on Gelfand pairs \cite{W2015b}. 
We sketch some of these developments.  Due to
constraints of time and space we pass over many aspects of operator theory
and orbit geometry, for example those described in \cite{BB2014},
\cite{BL2014} and \cite{BB2015}, related to
stepwise square integrable representations.
\medskip

In Section \ref{sec2} we recall the basic facts \cite{MW1973},
with a few extensions, on square integrable representations of nilpotent 
Lie groups.  
In Section \ref{sec3} we recall the concept and main results for
stepwise square integrable nilpotent Lie group.  
\medskip

In Section \ref{sec4} we show how nilradicals of minimal parabolic subgroups 
have the required decomposition for stepwise square integrability.  This 
is a construction based on concept of strongly orthogonal restricted roots.
\medskip

In Section \ref{sec5} we indicate the consequences for homogeneous compact 
nilmanifolds, and in Section \ref{sec6} we mention the application to 
analysis on commutative nilmanifolds.
\medskip

In Section \ref{sec7} we start the extension of stepwise square integrability
results from the nilradical $N$ of a minimal parabolic $P = MAN$ to various 
subgroups that contain $N$.  This section concentrates on the subgroup $MN$
and takes advantage of principal orbit theory.  That gives a sharp
simplification to the Plancherel and Fourier Inversion formulae.
In Section \ref{sec8} we look at $P$ and its subgroup $AN$.  They are 
not unimodular,
so we introduce the Dixmier-P\' ukanszky operator $D$ whose semi--invariance
balances that of the modular function.  It is a key point for the
Plancherel and Fourier Inversion formulae.
\medskip

Sections \ref{sec9} and \ref{sec10} are a short discussion of work in 
progress on the
extension of results from minimal parabolics to parabolics in general.
There are two places where matters diverge from the minimal parabolic
case.  First, there is a technical adjustment to the definition of stepwise
square integrable representation, caused by the fact that in the non--minimal
case the restricted roots need not form a root system.  Second, again for
technical reasons, the explicit Plancherel Formula only comes through
for the maximal amenable subgroups $UAN$ of $G$, and not for all of the
parabolic.
\medskip

This work was partially supported by a Simons Foundation grant and by
the award of a Dickson Emeriti Professorship.  It expands a talk at the
$11$-th International Workshop ``Lie Theory and Its Applications in Physics''
in Varna.  My thanks to Prof. Vladimir Dobrev and the others on the 
organizing committee for hospitality at that Workshop.

\section{\hskip -4pt .\hskip 2pt  Square Integrable Representations.}
\label{sec2}
\setcounter{lemma}{0}
\setcounter{theorem}{0}
\setcounter{proposition}{0}
\setcounter{corollary}{0}
\setcounter{definition}{0}
\setcounter{remark}{0}
\setcounter{example}{0}

Let $G$ be a unimodular locally compact group with center $Z$, and
let $\pi$ be an irreducible unitary representation. We associate the
central character $\chi_\pi \in \widehat{Z}$ by $\pi(z) = \chi_\pi(x)\cdot 1$
for $z \in Z$.  Consider a matrix coefficient $f_{u,v}: x \mapsto
\langle u, \pi(x)v \rangle$.  Then $|f_{u,v}|$ is a well defined function
on $G/Z$.  Fix Haar measures $\mu_G$ on $G$, $\mu_Z$ on $Z$ and
$\mu_{G/Z}$ on $G/Z$ such that $d\mu_G = d\mu_Z \, d\mu_{G/Z}$\,.
The following results are well known.

\begin{theorem}\label{sq}
The following conditions on $\pi \in \widehat{G}$ are equivalent.

{\rm (1)} There exist nonzero $u,v \in \cH_\pi$ with $|f_{u,v}| \in \cL^2(G/Z)$.

{\rm (2)} $|f_{u,v}| \in \cL^2(G/Z)$ for all $u,v \in \cH_\pi$.

{\rm (3)} 
$\pi$ is a discrete summand of the representation $\Ind_Z^G(\chi_\pi)$.
\end{theorem}

\begin{theorem} If the conditions of {\rm Theorem \ref{sq}} are satisfied for
an irreducible $\pi \in \widehat{G}$, then there is a number $\deg \pi > 0$
such that 
\begin{equation}
\int_{G/Z} f_{u,v}(x)\overline{f_{u',v'}(x)}d\mu_{G/Z}(xZ)
	= \tfrac{1}{\deg \pi} \langle u,u' \rangle
		\overline{\langle vv'\rangle}
\end{equation}
for all $u,u',v,v' \in \cH_\pi$\,.
If $\pi_1, \pi_2 \in \widehat{G}$ are inequivalent and satisfy the
conditions of {\rm Theorem \ref{sq}}, and $\chi_{\pi_1} = \chi_{\pi_2}$,
then
\begin{equation}
\int_{G/Z} \langle u, \pi_1(x)v \rangle 
	\overline{\langle u', \pi_2(x)v' \rangle}d\mu_{G/Z}(xZ) = 0
\end{equation}
for all $u,v \in \cH_{\pi_1}$ and all $u', v' \in \cH_{\pi_2}$\,.
\end{theorem}

The main results of \cite{MW1973} shows exactly how this works for 
nilpotent Lie groups:

\begin{theorem}  \label{mw-nilp}
Let $N$ be a connected simply connected Lie group with center $Z$,
$\gn$ and $\gz$ their Lie algebras, and $\gn^*$ the linear dual space 
of $\gn$.  Let $\lambda \in \gn^*$ and let $\pi_\lambda$ denote the
irreducible unitary representation attached to $\Ad^*(N)\lambda$ by
the Kirillov theory {\rm \cite{K1962}}.  Then the following conditions are 
equivalent.

{\rm (1)} $\pi_\lambda$ satisfies the conditions of {\rm Theorem \ref{sq}}.

{\rm (2)} The coadjoint orbit $\Ad^*(N)\lambda = \{\nu \in \gn^* \mid 
	\nu|_\gz = \lambda|_\gz$.

{\rm (3)} The bilinear form $b_\lambda(x,y) = \lambda([x,y])$ on $\gn/\gz$
	is nondegenerate.

{\rm (4)} The universal enveloping algebra $\cU(\gz)$ is the center of $\cU(\gn)$.

\noindent
The Pfaffian polynomial $\Pf(b_\lambda)$ is a polynomial function 
$P(\lambda|_{\gz})$ on $\gz^*$, and the set of representations 
$\pi_\lambda$ for which these conditions hold, is parameterized by the set
$\{\gamma \in \gz^* \mid P(\gamma) \ne 0\}$ $($which is empty or Zariski open 
in $\gz^*)$.
\end{theorem}

We will say that the connected simply connected Lie group $N$ is
{\em square integrable} if there exists $\lambda \in \gn^*$ such that
$P(\lambda|_{\gz}) \ne 0\}$.  For convenience we will sometimes write
$P(\lambda)$ for $P(\lambda|_{\gz})$ and $\pi_\gamma$ for $\pi_\lambda$
where $\gamma = \lambda|_{\gz}$\,.

\begin{theorem} \label{planch-conc}
Let $N$ be a square integrable connected simply connected Lie group
with center $Z$.  Then Plancherel measure on $\widehat{N}$ is 
concentrated on $\{\pi_\lambda \mid P(\lambda) \ne 0\}$, and there
the Plancherel measure is given by the measure $|P(\lambda)d\lambda|$
on $\gz^*$ and the formal degree $\deg \pi_\lambda = |P(\lambda|_\gz)|$.
\end{theorem}

Given $\gamma \in \gz^*$ with $P(\gamma) \ne 0$ and a Schwartz class
($\cC(N)$) function $f$ on $N$ we write $\cO(\gamma)$ for the co-adjoint orbit
$\Ad^*(N)\gamma = \gamma + \gz^\perp$\,, $f_\gamma$ for the restriction 
of $f\cdot \exp$ to $\cO(\gamma)$,
and $\widehat{f_\gamma}$ for the 
Fourier transform of $f_\gamma$ on $\cO(\gamma)$.

\begin{theorem} \label{sq-inv}
Let $N$ be a square integrable connected simply connected Lie group
with center $Z$ and $f \in \cC(N)$.  
If $\gamma \in \gz^*$ with $P(\gamma) \ne 0$ then the distribution
character of $\pi_\gamma$ is given by
\begin{equation}\label{char-sq}
\Theta_{\pi_\gamma}(f) = 
	\trace \int_N f(x) \pi_\gamma(x) d\mu_G(x) =
	c^{-1}|P(\gamma)|^{-1}\int_{\nu \in \cO(\gamma)} \widehat{f_\gamma} \,d\nu
\end{equation}
where $c = d!2^d \text{ and } d = \dim(\gn / \gz)/2$ and $d\nu$ is
ordinary Lebesgue measure on the affine space $\cO(\gamma)$\,.
The Fourier Inversion formula for $N$ is
\begin{equation}\label{fi-sq}
f(x) = c\int_{\gz^*} \Theta_\gamma(r_xf)
		|P(\gamma)|\, d\gamma \text{ where }
	(r_xf)(y) = f(yx) \text{ (right translate)}.
\end{equation}
\end{theorem}

There also are multiplicity results on $\cL^2(N/\Gamma)$ where 
$N$ is square integrable and
$\Gamma$ is a discrete co-compact subgroup, but they are the
same as in the stepwise square integrable case, so we postpone their
description.

\section{\hskip -4pt . Stepwise Square Integrability.}
\setcounter{lemma}{0}
\setcounter{theorem}{0}
\setcounter{proposition}{0}
\setcounter{corollary}{0}
\setcounter{definition}{0}
\setcounter{remark}{0}
\setcounter{example}{0}
\label{sec3}

In order to to go beyond square integrable nilpotent
groups, we suppose that the connected simply connected nilpotent Lie group
decomposes as
\begin{equation}\label{setup}
\begin{aligned}
N = &L_1L_2\dots L_{m-1}L_m \text{ where }\\
 &\text{(a) each $L_r$ has unitary representations with coeff in
$\cL^2(L_r/Z_r)$,} \\
 &\text{(b) } N_r := L_1L_2\dots L_r \text{ is normal in } N
   \text{ with } N_r = N_{r-1}\rtimes L_r\,, \\
 &\text{(c) } [\gl_r,\gz_s] = 0 \text{ and } [\gl_r,\gl_s] \subset \gv
        \text{ for } r > s \text{ with } \gl_r = \gz_r + \gv_r \\
 &\phantom{\text{(c) }}\text{ where } \gn = \gs + \gv,\,\, 
	\gs = \oplus \, \gz_r\, \text{ and } \gv = \oplus \, \gv_r\,.
\end{aligned}
\end{equation}
We will use the following notation.
\begin{equation}\label{c-d}
\begin{aligned}
&\text{(a) }d_r = \tfrac{1}{2}\dim(\gl_r/\gz_r) \text{ so }
        \tfrac{1}{2} \dim(\gn/\gs) = d_1 + \dots + d_m\,,\\
        & \phantom{XXXXXXX} \text{ and } 
		c = 2^{d_1 + \dots + d_m} d_1! d_2! \dots d_m!\\
&\text{(b) }b_{\lambda_r}: (x,y) \mapsto \lambda([x,y])
        \text{ viewed as a bilinear form on } \gl_r/\gz_r \\
&\text{(c) }S = Z_1Z_2\dots Z_m = Z_1 \times \dots \times Z_m \text{ where } Z_r
        \text{ is the center of } L_r \\
&\text{(d) }P: \text{ polynomial } P(\lambda) = \Pf(b_{\lambda_1})
        \Pf(b_{\lambda_2})\dots \Pf(b_{\lambda_m}) \text{ on } \gs^* \\
&\text{(e) }\gt^* = \{\lambda \in \gs^* \mid P(\lambda) \ne 0\} \\
&\text{(f) } \pi_\lambda \in \widehat{N} \text{ for } \lambda \in \gs^*
    \text{ with $P(\lambda) \ne 0$, irreducible unitary representation }\\
	&\phantom{XXXXXXX} \text{of } N = L_1L_2\dots L_m
        \text{ constructed as follows. }
\end{aligned}
\end{equation}

Start with the representation $\pi_{\lambda_1} \in \widehat{N_1}$ specified by 
$\lambda_1 \in \gz_1^*$ with $\Pf(b_{\lambda_1}) \ne 0$.  
Choose an invariant polarization $\gp_1' \subset \gn_2$ for the linear
functional $\lambda_1' \in \gn_2^*$ that agrees with $\lambda_1$ on $\gn_1$
and vanishes on $\gl_2$.  Since $L_r$ centralizes $S_{r-1}$,
$\ad^*(\gl_2)(\lambda_1')|_{\gz_1 + \gl_2} = 0$,
so  $\gp_1' = \gp_1 + \gl_2$ where $\gp_1$ is an
invariant polarization for $\lambda_1 \in \gn_1^*$.
The associated representations
are $\pi'_{\lambda_1} \in \widehat{N_2}$ and $\pi_{\lambda_1} \in
\widehat{N_1}$.  Note that $N_2/P_1' = N_1/P_1$\,,
so the representation spaces $\cH_{\pi'_{\lambda_1}} = \cL^2(N_2/P_1')
= \cL^2(N_1/P_1) = \cH_{\pi_{\lambda_1}}$.
In other words, $\pi'_{\lambda_1}$ extends $\pi_{\lambda_1}$ to a
unitary representation of $N_2$ on the same Hilbert space
$\cH_{\pi_{\lambda_1}}$, and $d\pi_{\lambda_1'}(\gz_2) = 0$.
Now the Mackey Little Group method gives us

\addtocounter{theorem}{1}
\begin{lemma}\label{ext-1-2}
The irreducible unitary representations of $N_2$, whose restrictions to
$N_1$ are multiples of $\pi_{\lambda_1}$, are the
$\pi'_{\lambda_1} \widehat\otimes\gamma$ where $\gamma \in \widehat{L_2}
= \widehat{N_2/N_1}$\,.
\end{lemma}

Given $\lambda_2 \in \gz_2^*$ with $\Pf(b_{\lambda_2}) \ne 0$ we have
$\pi_{\lambda_2} \in \widehat{L_2}$ with coefficients in $\cL^2(L_2/Z_2)$.
In the notation of Lemma \ref{ext-1-2} we define
\begin{equation}
\pi_{\lambda_1 + \lambda_2} \in \widehat{N_2} \text{ by }
\pi_{\lambda_1 + \lambda_2}
= \pi'_{\lambda_1} \widehat\otimes \pi_{\lambda_2}\,.
\end{equation}

\addtocounter{theorem}{1}
\begin{proposition}\label{L2.1}
The coefficients $f_{z,w}(xy)= \langle z,
\pi_{\lambda_1 + \lambda_2}(xy)w\rangle$ of $\pi_{\lambda_1 + \lambda_2}$
belong to $\cL^2(N_2/S_2)$, in fact satisfy
$
||f_{z,w}||^2_{\cL^2(N_r/S_r)} = \tfrac{||z||^2 ||w||^2}{\deg(\pi_{\lambda_1})
        \dots \deg(\pi_{\lambda_r})}\,.
$
\end{proposition}

Proposition \ref{L2.1} starts a recursion using
$N_{r+1} = N_r \rtimes L_{r+1}$.
We fix nonzero $\lambda_i \in \gz_i^*$ for $1 \leqq i \leqq r+1$,
and we start with the representation $\pi_{\lambda_1 + \dots + \lambda_r}$
constructed step by step from the square integrable representations
$\pi_{\lambda_i} \in \widehat{L_i}$ for $1 \leqq i \leqq r$.  The
representation space $\cH_{\pi_{\lambda_1 + \dots + \lambda_r}} =
\cH_{\pi_{\lambda_1}} \widehat\otimes \dots \widehat\otimes \cH_{\pi_{\lambda_r}}$.  The
coefficients of $\pi_{\lambda_1 + \dots + \lambda_r}$ have absolute
value in $\cL^2(N_r/S_r)$.  They satisfy
\begin{equation}
||f_{z,w}||^2_{\cL^2(N_r/S_r)} = \tfrac{||z||^2 ||w||^2}{\deg(\pi_{\lambda_1})
        \dots \deg(\pi_{\lambda_r})}\,.
\end{equation}
Then $\pi_{\lambda_1 + \dots + \lambda_r}$ extends to a representation
$\pi'_{\lambda_1 + \dots + \lambda_r}$ of $L_{r+1}$ on the same Hilbert
space $\cH_{\pi_{\lambda_1 + \dots + \lambda_r}}$, and it satisfies
$d\pi'_{\lambda_1 + \dots + \lambda_r}(\gz_{r+1}) = 0$.  As in
Lemma \ref{ext-1-2},

\addtocounter{theorem}{1}
\begin{lemma}\label{ext-r}
The irreducibles $\pi \in \widehat{N_{r+1}}$, whose restrictions to
$N_r$ are multiples of $\pi_{\lambda_1 + \dots + \lambda_r}$, are the
$\pi'_{\lambda_1 + \dots + \lambda_r} \widehat\otimes\gamma$ where
$\gamma \in \widehat{L_{r+1}} = \widehat{N_{r+1}/N_r}$\,.
\end{lemma}

As in Proposition \ref{L2.1}, define
$\pi_{\lambda_1 + \dots + \lambda_{r+1}} = \pi'_{\lambda_1 + \dots + \lambda_r}
\widehat\otimes \pi_{\lambda_{r+1}}$\,.  Then

\addtocounter{theorem}{1}
\begin{proposition}\label{L2.2}
The coefficients $f_{z,w}(x_1 \dots x_{r+1})= \langle z,
\pi_{\lambda_1 + \dots + \lambda_{r+1}}(x_1x_2 \cdots x_{r+1})w\rangle$ of
$\pi_{\lambda_1 + \dots + \lambda_{r+1}}$
belong to $\cL^2(N_{r+1}/S_{r+1})$, in fact satisfy
$$
||f_{z,w}||^2_{\cL^2(N_{r+1}/S_{r+1})} = \tfrac{||z||^2 ||w||^2}
{\deg(\pi_{\lambda_1}) \dots \deg(\pi_{\lambda_{r+1}})}\,.
$$
\end{proposition}

Since $\deg \pi_{\lambda_r} = |\Pf(b_{\lambda_r})|$,
Proposition \ref{L2.2} is the recursion step for our construction. Passing
to the end case $r+1 = m$ we see that Plancherel 
measure is concentrated on $\{\pi_\lambda \mid \lambda \in \gt^*\}$.
Using (\ref{setup})(c) to see that conjugation by elements of $L_s$ 
has no effect on the $\Pf(b_{\lambda_r})$ for $r < s$, we arrive at

\begin{theorem}\label{plancherel-general}
Let $N$ be a connected simply connected nilpotent Lie group that
satisfies {\rm (\ref{setup})}.  Then Plancherel measure for $N$ is
concentrated on $\{\pi_\lambda \mid \lambda \in \gt^*, P(\lambda) \ne 0\}$.
If $\lambda \in \gt^*$, $P(\lambda) \ne 0$ 
and $u, v \in \cH_{\pi_\lambda}$\,, then
the coefficient $f_{u,v}(x) = \langle u, \pi_\nu(x)v\rangle$
satisfies
\begin{equation}
||f_{u,v}||^2_{\cL^2(N / S)} = {||u||^2||v||^2}/{|P(\lambda)|}\,.
\end{equation}
The distribution character 
$\Theta_{\pi_\lambda} : f \mapsto \tr \int_G f(x)\pi(x)dx$ of $\pi_{\lambda}$ 
is given by
\begin{equation}
\Theta_{\pi_\lambda}(f) = c^{-1}|P(\lambda)|^{-1}\int_{\cO(\lambda)}
        \widehat{f_\lambda}(\xi)d\nu_\lambda(\xi) \text{ for } f \in \cC(N)
\end{equation}
where $\cC(N)$ is the Schwartz space, 
$\cO(\lambda) = \Ad^*(N)\lambda = \gs^\perp + \lambda$,
$f_\lambda$ is the lift $f_\lambda(\xi) = f(\exp(\xi))$, 
$\widehat{f_\lambda}$ is its classical Fourier transform,
and $d\nu_\lambda$ is the translate of normalized Lebesgue measure from
$\gs^\perp$ to $\Ad^*(N)\lambda$.  Further,
\begin{equation}
f(x) = c\int_{\gt^*} \Theta_{\pi_\lambda}(r_xf) |P(\lambda)|d\lambda
        \text{ for } f \in \cC(N).
\end{equation}
\end{theorem}
\addtocounter{theorem}{1}
\begin{definition}\label{stepwise2}
{\rm The representations $\pi_\lambda$ of (\ref{c-d}(f)) are the
{\it stepwise square integrable} representations of $N$ relative to
(\ref{setup}).}\hfill $\diamondsuit$
\end{definition}

The left action $(l(x)f)(g) = f(x^{-1}g)$ and the right action $(r(y)f)(g)
= f(gy)$ of $N$ on functions carries over to coefficients of $\pi$ as
$l(x)r(y)f_{u,v} = f_{\pi(x)u,\pi(y)v}$.  If $\pi = \pi_\lambda$
stepwise square integrable, $u, v \in \cH_{\pi_\lambda}$ are $C^\infty$
vectors, and if $\Phi$ and $\Psi$ belong to the universal enveloping
algebra $\cU(\gn)$, then $l(\Phi)r(\Psi)f_{u,v} = f_{d\pi(\Psi)u,d\pi(\Phi)v}$
is just another coefficient, $C^\infty$ and $\cL^2(N/S)$.  If $\zeta_\lambda
\in \widehat{S}$ is the quasicentral character of $\pi_\lambda$ it follows 
that $f_{u,v}$ belongs to the relative Schwartz space $\cC(N/S,\zeta_\lambda)$.
In particular it follows that $|f_{u,v}| \in \cL^p(N/S)$ for all $p \geqq 1$.
Taking Schwartz class wave packets over $S$ of coefficient functions of
stepwise square integrable representations of $N$ one can express the
Plancherel formula of Theorem \ref{plancherel-general} in terms of
coefficient functions.

\section{\hskip -4pt .\hskip 2pt  Nilradicals of Minimal Parabolics.}
\setcounter{lemma}{0}
\setcounter{theorem}{0}
\setcounter{proposition}{0}
\setcounter{corollary}{0}
\setcounter{definition}{0}
\setcounter{remark}{0}
\setcounter{example}{0}
\label{sec4}
Fix a real simple Lie group $G$, an Iwasawa decomposition $G = KAN$,
and a minimal parabolic subgroup $Q = MAN$ in $G$.
Let $m = \rank_\R G = \dim_\R A$\,.  As usual, write $\gk$ for the Lie
algebra of $K$, $\ga$ for the Lie algebra of $A$, and $\gn$ for the
Lie algebra of $N$.  Complete $\ga$ to a Cartan subalgebra $\gh$ of $\gg$.
Then $\gh = \gt + \ga$ with $\gt = \gh \cap \gk$.  Now we have root systems
\begin{equation}\label{root-systems}
\begin{aligned}
&\Delta(\gg_\C,\gh_\C)\text{: roots of } \gg_\C \text{ relative to } \gh_\C
	\text{ (ordinary roots),}\\
&\Delta(\gg,\ga)\text{: roots of } \gg \text{ relative to } \ga 
	\text{ (restricted roots),}\\
&\Delta_0(\gg,\ga) = \{\alpha \in \Delta(\gg,\ga) \mid
        2\alpha \notin \Delta(\gg,\ga)\} \text{ (nonmultipliable)}.
\end{aligned}
\end{equation}
Here $\Delta(\gg,\ga)$ and $\Delta_0(\gg,\ga)$ are root
systems in the usual sense.  Any positive root system
$\Delta^+(\gg_\C,\gh_\C) \subset \Delta(\gg_\C,\gh_\C)$ defines positive
systems
\begin{equation}
\begin{aligned}
&\Delta^+(\gg,\ga) = \{\gamma|_\ga \mid \gamma \in
	\Delta^+(\gg_\C,\gh_\C) \text{ and } \gamma|_\ga \ne 0\},\\
&\Delta_0^+(\gg,\ga) = \Delta_0(\gg,\ga) \cap \Delta^+(\gg,\ga).
\end{aligned}
\end{equation}
We can (and do) choose $\Delta^+(\gg,\gh)$ so that
\begin{equation}
\begin{aligned}
&\gn \text{ is the sum of the positive restricted root spaces and}\\
&\text{if } \gamma \in \Delta(\gg_\C,\gh_\C) \text{ and } \gamma|_\ga \in
	\Delta^+(\gg,\ga) \text{ then } \gamma \in \Delta^+(\gg_\C,\gh_\C).
\end{aligned}
\end{equation}

Two roots are called {\em strongly orthogonal} if their sum and their
difference are not roots.  Then they are orthogonal.  The Kostant cascade
construction is
\begin{equation}
\begin{aligned}
&\beta_1 \in \Delta^+(\gg,\ga) \text{ is a maximal positive restricted root
and }\\
& \beta_{r+1} \in \Delta^+(\gg,\ga) \text{ is a maximum among the roots of }
\Delta^+(\gg,\ga) \\ 
& \text{ that are orthogonal to all } \beta_i \text{ with } i \leqq r
\end{aligned}
\end{equation}
Then the $\beta_r$ are mutually strongly orthogonal.
Each $\beta_r \in \Delta_0^+(\gg,\ga)$, and
$\beta_1$ is unique because $\Delta(\gg,\ga)$ is irreducible.
For $1\leqq r \leqq m$ define
\begin{equation}\label{layers}
\begin{aligned}
&\Delta^+_1 = \{\alpha \in \Delta^+(\gg,\ga) \mid \beta_1 - \alpha \in \Delta^+(\gg,\ga)\}
\text{ and }\\
&\Delta^+_{r+1} = \{\alpha \in \Delta^+(\gg,\ga) \setminus (\Delta^+_1 \cup \dots \cup \Delta^+_r)
        \mid \beta_{r+1} - \alpha \in \Delta^+(\gg,\ga)\}.
\end{aligned}
\end{equation}
\addtocounter{theorem}{1}
\begin{lemma} \label{fill-out}
If $\alpha \in \Delta^+(\gg,\ga)$, either
$\alpha \in \{\beta_1, \dots , \beta_m\}$
or $\alpha$ belongs to just one $\Delta^+_r$\,.
\end{lemma}
\addtocounter{theorem}{1}
\begin{lemma}\label{layers2}
$\Delta^+_r\cup \{\beta_r\}
= \{\alpha \in \Delta^+ \mid \alpha \perp \beta_i \text{ for } i < r
\text{ and } \langle \alpha, \beta_r\rangle > 0\}.$
In particular, $[\gl_r,\gl_s] \subset \gl_t$ where $t = \min\{r,s\}$.
\end{lemma}
Lemma \ref{fill-out} shows that the Lie algebra $\gn$ of $N$ is the
direct sum of its subspaces
\begin{equation}\label{def-m}
\gl_r = \gg_{\beta_r} + {\sum}_{\Delta^+_r}\, \gg_\alpha
\text{ for } 1\leqq r\leqq m
\end{equation}
and Lemma \ref{layers2} shows that $\gn$ has an increasing foliation
by ideals
\begin{equation}\label{def-filtration}
\gn_r = \gl_1 + \gl_2 + \dots + \gl_r \text{ for } 1 \leqq r \leqq m.
\end{equation}
Now we will see that the corresponding group level decomposition
$N = L_1L_2\dots L_m$ and the semidirect product decompositions
$N_r = N_{r-1}\rtimes L_r$ satisfy (\ref{setup}).  Denote
\begin{equation}
\begin{aligned}
&s_{\beta_r} \text{ is the Weyl group reflection in } \beta_r
\text{ and } \\
&\sigma_r: \Delta(\gg,\ga) \to \Delta(\gg,\ga) \text{ by }
\sigma_r(\alpha) = -s_{\beta_r}(\alpha).
\end{aligned}
\end{equation}
Note that $\sigma_r(\beta_s) = -\beta_s$ for $s \ne r$, $+\beta_s$ if $s = r$.
If $\alpha \in \Delta^+_r$ we still have $\sigma_r(\alpha) \perp \beta_i$
for $i < r$ and $\langle \sigma_r(\alpha), \beta_r\rangle > 0$.  If
$\sigma_r(\alpha) < 0$ then $\beta_r - \sigma_r(\alpha) > \beta_r$
contradicting maximality of $\beta_r$.  Thus, using
Lemma \ref{layers2}, $\sigma_r(\Delta^+_r) = \Delta^+_r$.

\addtocounter{theorem}{1}
\begin{lemma} \label{layers-nilpotent}
If $\alpha \in \Delta^+_r$ then $\alpha + \sigma_r(\alpha) = \beta_r$.
{\rm (}It is possible that
$\alpha = \sigma_r(\alpha) = \tfrac{1}{2}\beta_r$ when
$\tfrac{1}{2}\beta_r$ is a root.{\rm ).}
If $\alpha, \alpha' \in \Delta^+_r$ and $\alpha + \alpha' \in \Delta(\gg,\ga)$
then $\alpha + \alpha' = \beta_r$\,.
\end{lemma}
\addtocounter{theorem}{1}
\begin{lemma}\label{pairing}
Let $\gn$ be a nilpotent Lie algebra, $\gz$ its center, and $\gv$ a
vector space complement to $\gz$ in $\gn$.  Suppose that 
$\gv = \gu + \gu'$, $\gu = \sum \gu_a$ and $\gu' = \sum \gu'_a$\,, and
$\gz = \sum \gz_b$ with $\dim \gz_b = 1$ in such a way that
{\rm (i)} each $[\gu_a, \gu_a] = 0 = [\gu'_a,\gu'_a]$,
{\rm (ii)} if $a_1 \ne a_2$ then $[\gu_{a_1} , \gu'_{a_2}] = 0$ and
{\rm (iii)} for each $a$ there is a nondegenerate pairing 
$\gu_a \otimes \gu'_a \to \gz_{b_a}$\,, by $u \otimes u' \mapsto [u,u']$.
Then $\gn$ is a direct sum of Heisenberg
algebras  $\gz_{b_a} + \gu_a + \gu'_a$ and the commutative algebra that is the
sum of the remaining $\gz_b$\,.
\end{lemma}

Now one runs through a number of special situations:
(1) If $\gg$ is the split real form of $\gg_\C$ then each $L_r$ has square
integrable representations.
(2) If $\gg$ is simple but not absolutely simple then each $L_r$ has square
integrable representations.
(3) If $G$ is the quaternion special linear group $SL(n;\H)$ then
$L_1$ has square integrable representations.
(4) If $G$ is the group $E_{6,F_4}$ of collineations of the Cayley
projective plane then $L_1$ has square integrable representations.
(5) The group $L_1$ has square integrable representations.
(6) If $\gg$ is absolutely simple then each $L_r$ has square
integrable representations.
Putting these together, Theorem \ref{plancherel-general}
applies to nilradicals of minimal parabolic subgroups:

\begin{theorem}\label{iwasawa-layers}
Let $G$ be a real reductive Lie group, $G = KAN$ an Iwasawa
decomposition, $\gl_r$ and $\gn_r$ the subalgebras of $\gn$ defined in
{\rm (\ref{def-m})} and {\rm (\ref{def-filtration})},
and $L_r$ and $N_r$ the corresponding analytic subgroups of $N$.
Then the $L_r$ and $N_r$ satisfy {\rm (\ref{setup})}.  In particular,
Plancherel measure for $N$ is
concentrated on $\{\pi_\lambda \mid \lambda \in \gt^*\}$.
If $\lambda \in \gt^*$, and if $u$ and $v$ belong to the
representation space $\cH_{\pi_\lambda}$ of $\pi_\lambda$,  then
the coefficient $f_{u,v}(x) = \langle u, \pi_\lambda(x)v\rangle$
satisfies
$
||f_{u,v}||^2_{\cL^2(N / S)} = \frac{||u||^2||v||^2}{|P(\lambda)|}\,.
$
The distribution character $\Theta_{\pi_\lambda}$ of $\pi_{\lambda}$ satisfies
$
\Theta_{\pi_\lambda}(f) = c^{-1}|P(\lambda)|^{-1}\int_{\cO(\lambda)}
        \widehat{f_\lambda}(\xi)d\nu_\lambda(\xi) \text{ for } f \in \cC(N).
$
Here $\cC(N)$ is the Schwartz space, 
$\cO(\lambda)$ is the coadjoint orbit $\Ad^*(N)\lambda = \gs^\perp + \lambda$,
$f_\gamma$ is the lift $f_\gamma(\xi) = f(\exp(\xi))$ to $\gs^\perp +\lambda$, 
$\widehat{f_\gamma}$ is its classical Fourier transform ,
and $d\nu_\lambda$ is the translate of normalized Lebesgue measure from
$\gs^\perp$ to $\Ad^*(N)\lambda$.  The Plancherel formula on $N$ is
$
f(x) = c\int_{\gt^*} \Theta_{\pi_\lambda}(r_xf) |P(\lambda)|d\lambda
        \text{ for } f \in \cC(N).
$
\end{theorem}

\section{\hskip -4pt .\hskip 2pt  Compact Nilmanifolds.}
\setcounter{lemma}{0}
\setcounter{theorem}{0}
\setcounter{proposition}{0}
\setcounter{corollary}{0}
\setcounter{definition}{0}
\setcounter{remark}{0}
\setcounter{example}{0}
\label{sec5}

Here are the basic facts on discrete uniform (i.e. co-compact) subgroups
of connected simply connected nilpotent Lie groups.
See \cite[Chapter 2]{R1972} for an exposition.

\addtocounter{theorem}{1}
\begin{proposition}\label{alg-1}
The following are equivalent.
\begin{itemize}
\item $N$ has a discrete subgroup $\Gamma$ with $N/\Gamma$ compact.
\item $N \cong N_\R$ where $N_\R$ is the group of real points in a unipotent
  linear algebraic group defined over the rational number field $\Q$
\item $\gn$ has a basis $\{\xi_j\}$ for which
the coefficients $c_{i,j}^k$ in
$[\xi_i,\xi_j] = \sum c_{i,j}^k\xi_k$ are rational numbers.
\end{itemize}
Under those conditions let $\gn_{_\Q}$ denote the rational span of
$\{\xi_j\}$ and let $\gn_{_\Z}$ be the integral span.  Then
$\exp(\gn_{_\Z})$ generates a discrete subgroup $N_\Z$ of $N = N_\R$
and $N_\R/N_\Z$ is compact.  Conversely, if $\Gamma$ is a discrete
co-compact subgroup of $N$ then the $\Z$--span of $\exp^{-1}(\Gamma)$
is a lattice in $\gn$ for which any generating set $\{\xi_j\}$ is a
basis of $\gn$ such that the coefficients $c_{i,j}^k$ in
$[\xi_i,\xi_j] = \sum c_{i,j}^k\xi_k$ are rational numbers.
\end{proposition}

The conditions of Proposition \ref{alg-1} hold for the nilpotent
groups
studied in Section \ref{sec4}; there one can choose the basis
$\{\xi_j\}$ of $\gn$ so that the $c_{i,j}^k$ are integers.
\smallskip

The basic facts on square integrable representations that occur in
compact quotients $N/\Gamma$,
as described in \cite[Theorem 7]{MW1973}, are

\addtocounter{theorem}{1}
\begin{proposition}\label{mult-1}
Let $N$ be a connected simply connected nilpotent Lie group that has
square integrable representations, and let $\Gamma$ a discrete co-compact
subgroup.  Let $Z$ be the center of $N$ and normalize the volume form on
$\gn / \gz$ by normalizing Haar measure on $N$ so that $N/Z\Gamma$ has
volume $1$.  Let $P$ be the corresponding
Pfaffian polynomial on $\gz^*$.  Note that $\Gamma \cap Z$ is a lattice
in $Z$ and $\exp^{-1}(\Gamma \cap Z)$ is a lattice (denote it $\Lambda$)
in $\gz$.  That defines the dual lattice $\Lambda^*$ in $\gz^*$.  Then
a square integrable representation $\pi_\lambda$ occurs in $\cL^2(N/\Gamma)$
if and only if $\lambda \in \Lambda^*$, and in that case $\pi_\lambda$
occurs with multiplicity $|P(\lambda)|$.
\end{proposition}

\addtocounter{theorem}{1}
\begin{definition}\label{ref-rational}
{\rm
Let $N = N_\R$ be defined over $\Q$ as in Proposition \ref{alg-1}, so we
have a fixed rational form $N_\Q$.  We say that a connected Lie subgroup
$L \subset N$ is {\em rational} if $L \cap N_\Q$ is a rational form of
$L$, in other words if $\gl \cap \gn_{_\Q}$ contains a basis of $\gl$.
We say that a decomposition {\rm (\ref{setup})} is {\em rational} if the
subgroups $L_r$ and $N_r$ are rational.
}
\hfill $\diamondsuit$
\end{definition}

The following is immediate from this definition.
\addtocounter{theorem}{1}
\begin{lemma}\label{immediate}
Let $N$ be defined over $\Q$ as in {\rm Proposition \ref{alg-1}}
with rational structure defined by a discrete co-compact subgroup $\Gamma$.
If the decomposition {\rm (\ref{setup})} is rational then each
$\Gamma \cap Z_r$ in $Z_r$\,,
each $\Gamma \cap L_r$ in $L_r$\,,
each $\Gamma \cap S_r$ in $S_r$\,,
and each $\Gamma \cap N_r$ in $N_r$\,, is a
discrete co-compact subgroup defining the same rational structure as
the one defined by its intersection with $N_\Q$\,.
\end{lemma}

Now assume that $N$ and $\Gamma$ satisfy
the rationality conditions of Lemma \ref{immediate}.  Then for each $r$,
$Z_r \cap \Gamma$ is a lattice in the center $Z_r$ of $L_r$, and
$\Lambda_r := \log(Z_r \cap \Gamma)$ is a lattice in its Lie algebra $\gz_r$.
That defines the dual lattice $\Lambda_r^*$ in $\gz_r^*$.  We normalize the
Pfaffian polynomials on the $\gz_r^*$, and thus the polynomial $P$ on
$\gs^*$, by requiring that the $N_r/(S_r\cdot (N_r\cap\Gamma))$ have volume $1$.

\begin{theorem}\label{mult-form}
Let $\lambda \in \gt^*$.  Then a stepwise square integrable
representation $\pi_\lambda$ of $N$ occurs in $\cL^2(N/\Gamma)$ if and
only if each $\lambda_r \in \Lambda_r^*$\,, and
in that case the multiplicity of $\pi_\lambda$ on $\cL^2(N/\Gamma)$ is
$|P(\lambda)|$.
\end{theorem}

\section{\hskip -4pt .\hskip 2pt  Commutative Spaces.}
\setcounter{lemma}{0}
\setcounter{theorem}{0}
\setcounter{proposition}{0}
\setcounter{corollary}{0}
\setcounter{definition}{0}
\setcounter{remark}{0}
\setcounter{example}{0}
\label{sec6}

A commutative space $X = G/K$, or equivalently a Gelfand pair $(G,K)$,
consists of a locally compact group $G$ and a compact subgroup $K$ such
that the convolution algebra $\cL^1(K\backslash G/K)$ is commutative.
When $G$ is a connected Lie group it is equivalent to say that the algebra
$\cD(G,K)$ of $G$--invariant differential operators on $G/K$ is commutative.
We say that the commutative space $G/K$ is a commutative nilmanifold
if it is a nilmanifold in the sense that some nilpotent analytic
subgroup $N$ of $G$ acts transitively.  When $G/K$ is connected and
simply connected it follows that $N$ is the nilradical of $G$, that $N$
acts simply transitively on $G/K$, and that $G$ is the semidirect
product group $N\rtimes K$, so that $G/K = (N\rtimes K)/K$.  In this
section we study the role of square integrability and stepwise square 
integrability for commutative nilmanifolds $G/K = (N\rtimes K)/K$.
\medskip

The cases where $G/K$ and $(G,K)$ are
{\em irreducible} in the sense that $[\gn,\gn]$ (which must be central)
is the center of $\gn$ and $K$ acts irreducibly on $\gn/[\gn,\gn]$,
have been classified by E. B. Vinberg (\cite{V2001}, \cite{V2003}).
See \cite[\S 13.4B]{W2007} for the Lie algebra structure 
$\gv \times \gv \to \gz$.
The classification of commutative nilmanifolds is based on Vinberg's work
and was completed by O. Yakimova in \cite{Y2005} and 
\cite{Y2006}.
\medskip

It turns out that almost all commutative manifolds correspond to
nilpotent groups that are square integrable.  The exceptions are those
with a certain direct factor, and in those cases the nilpotent group is 
stepwise square integrable in two steps, so in those cases the
Plancherel formula follows directly from the general result above.
See \cite{W2015b} for the details.

\section{\hskip -4pt .\hskip 2pt  Minimal Parabolics: Subgroup $MN$.}
\label{sec7}
\setcounter{lemma}{0}
\setcounter{theorem}{0}
\setcounter{proposition}{0}
\setcounter{corollary}{0}
\setcounter{definition}{0}
\setcounter{remark}{0}
\setcounter{example}{0}
Fix an Iwasawa decomposition $G = KAN$ for a simple Lie group $G$ and the
minimal parabolic subgroup $Q = MAN$.  As usual, write $\gk$ for the Lie 
algebra of $K$, $\ga$ for the Lie algebra of $A$, $\gm$ for the Lie
algebra of $M$, and $\gn$ for the
Lie algebra of $N$.  Complete $\ga$ to a Cartan subalgebra $\gh$ of $\gg$.
Then we have root systems $\Delta(\gg_\C,\gh_\C)$, $\Delta(\gg,\ga)$ and 
$\Delta_0(\gg,\ga)$ described in (\ref{root-systems}).  $M$ is the
centralizer of $A$ in $K$.  Write ${}^0$ for identity component; 
then $Q^0 = M^0AN$.
\medskip

Recall the $\Pf$--nonsingular 
set $\gt^* = \{\lambda \in \gs^* \mid \Pf(b_\lambda) \ne 0\}$ of
{\rm (\ref{c-d}e)}; so $\Ad^*(M)\gt^* = \gt^*$.  Further, if 
$\lambda \in \gt^*$ and $c \ne 0$ then $c\lambda \in \gt^*$, 
in fact $\Pf(b_{c\lambda}) = c^{\dim(\gn/\gs)/2} \Pf(b_\lambda)$.
\medskip

Fix an $M$--invariant inner product $(\mu,\nu)$ on $\gs^*$\,. So $\Ad^*(M)$
preserves each sphere 
$\gs^*_t = \{\lambda \in \gs^* \mid (\lambda,\lambda) = t^2\}$.
Two orbits
$\Ad^*(M)\mu$ and $\Ad^*(M)\nu$ are of the {\sl same orbit type} if the 
isotropy subgroups $M_\mu$ and $M_\nu$ are conjugate, and an orbit is
{\sl principal} if all nearby orbits are of the same type.  Since $M$ and
$\gs^*_t$ are compact, there are only finitely many orbit types of $M$
on $\gs^*_t$, there is only one principal orbit type, and the union of the
principal orbits forms a dense open subset of $\gs^*_t$ whose complement has
codimension $\geqq 2$.  See \cite[Chapter 4, Section 3]{B1972} for a complete 
treatment of this material, or  
\cite[Part II, Chapter 3, Section 1]{G1993} modulo
references to \cite{B1972}, or \cite[Cap. 5]{N2005} for a basic
treatment, still with some references to \cite{B1972}.
\medskip

The action of $M$ on $\gs^*$ commutes with dilation so the
structural results on the $\gs_t$ also hold on 
$\gs^* = \bigcup_{t > 0} \gs^*_t$.  
Define the $\Pf$-nonsingular principal orbit set as follows:
\begin{equation}\label{defregset}
\gu^* = \{\lambda \in \gt^* \mid \Ad^*(M)\lambda \text{ is a principal }
        M\text{-orbit on } \gs^*\}.
\end{equation} 
Now principal orbit set $\gu^*$ is a dense open set with complement 
of codimension $\geqq 2$
in $\gs^*$.  If $\lambda \in \gu^*$ and $c \ne 0$ then
$c\lambda \in \gu^*$ with isotropy $M_{c\lambda} = M_\lambda$\,.
If $\lambda \in \gu^*_t := \gu^* \cap \gs^*_t$\,, so $\Ad^*(M)\lambda$ 
is a $\Pf$-nonsingular principal orbit of $M$ on the sphere $\gs^*_t$, then 
$\Ad^*(M^0)\lambda$ is a principal orbit of $M^0$ on $\gs^*_t$.  
Principal orbit isotropy subgroups of compact connected
linear groups are studied in \cite{HH1970} and the possibilities for 
the isotropy $(M^0)_\lambda$ are essentially known.  The following lets
us go from $(M^0)_\lambda$ to $M_\lambda$\,.

\addtocounter{theorem}{1}
\begin{proposition}\label{m-components-2} {\rm (\cite{W2014})}
Suppose that $G$ is connected and linear.  Then $M = F Z_G M^0$ where 
$Z_G$ is the center of $G$, $F = (\exp(i\ga) \cap K)$ is an elementary 
abelian $2$--group, and $Ad^*(F)$ acts trivially on $\gs^*$.  If
$\lambda \in \gu^*$ then the isotropy $M_\lambda = F Z_G (M^0)_\lambda$\,.
\end{proposition}

Thus the groups $M_\lambda$ are specified by the
work of W.--C. and W.--Y. Hsiang \cite{HH1970}.
\medskip

Given $\lambda \in \gu^*$ the stepwise square integrable representation
$\pi_\lambda \in \widehat{N}$ one proves that the Mackey 
obstruction $\varepsilon \in H^2(M_\lambda;U(1))$ is trivial, and
in fact that 
$\pi_\lambda$ extends 
to a unitary representation $\pi_\lambda^\dagger$ of $N\rtimes M_\lambda$ 
on the representation space of $\pi_\lambda$\,.
\medskip

Each $\lambda \in \gu^*$ now defines classes
\begin{equation}\label{nm-lambda-family}
 \cE(\lambda) := \left \{\pi_\lambda^\dagger \otimes \gamma \mid
\gamma \in \widehat{M_\lambda}\right \},\,\,
 \cF(\lambda) := \left \{\Ind_{NM_\lambda}^{NM} 
(\pi_\lambda^\dagger \otimes \gamma )\mid \pi_\lambda^\dagger
\otimes \gamma \in \cE(\lambda)\right \}
\end{equation} 
of irreducible unitary representations of $N\rtimes M_\lambda$ and
$NM$.  The Mackey little group method, plus the fact that the
Plancherel density on $\widehat{N}$ is polynomial on $\gs^*$\,, 
and $\gs^*\setminus\gu^*$ has measure $0$ in $\gt^*$, 
gives us

\addtocounter{theorem}{1}
\begin{proposition}\label{rep-mn}
Plancherel measure for $NM$ is concentrated on 
$\bigcup_{\lambda \in \gu^*}\cF(\lambda)$, equivalence classes
of irreducible representations 
$\eta_{\lambda,\gamma} := \Ind_{NM_\lambda}^{NM} 
(\pi_\lambda^\dagger \otimes \gamma)$ such that
$\pi_\lambda^\dagger \otimes \gamma \in \cE(\lambda)$ and 
$\lambda \in \gu^*$.  Further
$$
\eta_{\lambda,\gamma}|_N = 
\left . \left ( \Ind_{NM_\lambda}^{NM}(\pi_\lambda^\dagger 
        \otimes \gamma)\right ) \right |_N = \int_{M/M_\lambda} 
        (\dim \gamma)\, \pi_{\Ad^*(m)\lambda}\, d(mM_\lambda).
$$
\end{proposition}

There is a Borel section $\sigma$ to
$\gu^* \to \gu^*/\Ad^*(M)$ that picks out an element in each $M$-orbit
so that $M$ has the same isotropy subgroup at each of those elements.  In
other words in each $M$-orbit on $\gu^*$ we measurably choose an element
$\lambda = \sigma(\Ad^*(M)\lambda)$ such that those isotropy subgroups
$M_\lambda$ are all the same.  Let us denote
\begin{equation}\label{m-diamond}
M_\diamondsuit \text{: isotropy subgroup of } M \text{ at }
\sigma(\Ad^*(M)\lambda) \text{ for every } \lambda \in \gu^*
\end{equation}
We replace $M_\lambda$ by $M_\diamondsuit$, independent of $\lambda
\in \gu^*$, in Proposition \ref{rep-mn}.  That lets us assemble to 
representations of Proposition \ref{rep-mn} for a Plancherel Formula, as 
follows.
Since $M$ is compact, we have the Schwartz space $\cC(NM)$ just as in the
discussion of $\cC(N)$ between (\ref{c-d}) and Theorem \ref{plancherel-general},
except that the pullback $\exp^*\cC(NM) \ne \cC(\gn + \gm)$.  The
same applies to $\cC(NA)$ and $\cC(NAM)$

\addtocounter{theorem}{1}
\begin{proposition}\label{planch-mn}
Let $f \in \cC(NM)$ and write $(f_m)(n) = f(nm) = ({}_nf)(m)$ for $n \in N$
and $m \in M$.  The Plancherel density at $\Ind_{NM_\diamondsuit}^{NM} 
(\pi_\lambda^\dagger \otimes \gamma)$ is $(\dim \gamma)|\Pf(b_\lambda)|$ 
and the Plancherel Formula for $NM$ is
$$
f(nm) = c\int_{\gu^*/\Ad^*(M)}\, 
        \sum_{\cF(\lambda)}
        \tr \eta_{\lambda,\gamma}(_n f _m)\cdot 
        \dim(\gamma)\cdot |\Pf(b_\lambda)|d\lambda
$$
where $c = 2^{d_1 + \dots + d_m} d_1! d_2! \dots d_m!$\,, from 
{\rm (\ref{c-d})},
as in {\rm Theorem \ref{plancherel-general}}.
\end{proposition}

\section{\hskip -4pt .\hskip 2pt  Minimal Parabolics: $MAN$ and $AN$.}
\label{sec8}
\setcounter{lemma}{0}
\setcounter{theorem}{0}
\setcounter{proposition}{0}
\setcounter{corollary}{0}
\setcounter{definition}{0}
\setcounter{remark}{0}
\setcounter{example}{0}

Let $G$ be a separable locally compact group of type I.
Then \cite[\S 1]{LW1978} the Plancherel formula for $G$
has form
\begin{equation}\label{LW}
f(x) = \int_{\widehat{G}} \tr\pi(D(r(x)f)) d\mu_{_G}(\pi)
\end{equation}
where $D$ is an invertible positive self adjoint operator on $L^2(G)$,
conjugation--semi-invariant of weight equal to the modular function $\delta_G$,
and $\mu$ is a positive Borel measure on the unitary dual $\widehat{G}$.
If $G$ is unimodular then $D$ is the identity and (\ref{LW}) reduces to the
usual Plancherel formula.  The point is that semi-invariance of $D$
compensates any lack of unimodularity.  See \cite[\S 1]{LW1978} for
a detailed discussion.
$D \otimes \mu$ is unique (up to normalization of Haar measures) and
one tries to find a ``best'' choice of $D$.  Given any such pair $(D,\mu)$
we refer to $D$ as a {\sl Dixmier--Puk\' anszky Operator} on $G$ and
to $\mu$ as the associated {\sl Plancherel measure} on $\widehat{G}$.
We will construct a  Dixmier--Puk\' anszky Operator from the
Pfaffian polynomial $\Pf(b_\lambda)$.
\medskip

Let $\delta_{AN}$ and $\delta_Q$ denote the modular functions on $AN$ and
on $Q = MAN$.  As $M$ is compact and $\Ad_Q(N)$ is unipotent on $\gp$, 
they are determined by their restrictions to  $A$, where they are
given by $\delta(\exp(\xi)) = \exp(\tr(\ad(\xi)))$ with $\xi = \log a \in \ga$.

\addtocounter{theorem}{1}
\begin{lemma}\label{trace}
Let $\xi \in \ga$.  Then $\frac{1}{2}(\dim \gl_r + \dim \gz_r) \in \Z$
for $1\leqq r\leqq m$ and\\
\phantom{ii}{\rm (i)} the trace of $\ad(\xi)$ on $\gl_r$ is
  $\frac{1}{2}(\dim \gl_r + \dim \gz_r)\beta_r(\xi)$, \\
\phantom{i}{\rm (ii)} the trace of $\ad(\xi)$ on $\gn$ and on $\gp$ is
  $\frac{1}{2}\sum_r (\dim \gl_r + \dim \gz_r)\beta_r(\xi)$, \\
{\rm (iii)}
  the determinant of $\Ad(\exp(\xi))$ on $\gn$ and on $\gp$ is
  $\prod_r \exp(\beta_r(\xi))^{\frac{1}{2} (\dim \gl_r + \dim \gz_r)}$,\\
{\rm (iv)}
  $\delta_Q(man)
= \prod_r \exp(\beta_r(\log a))^{\frac{1}{2} (\dim \gl_r + \dim \gz_r)}$
and $\delta_{AN} = \delta_Q|_{AN}$.
\end{lemma}
Now compute

\addtocounter{theorem}{1}
\begin{lemma}\label{trace2}
Let $\xi \in \ga$ and $a = \exp(\xi) \in A$.  Then
$\ad(\xi)\Pf = \left (\frac{1}{2}\sum_r\, \dim (\gl_r/\gz_r) \beta_r(\xi)\right )\Pf$
and $\Ad(a)\Pf =
\left ( \prod_r \exp(\beta_r(\xi))^{\frac{1}{2}\dim (\gl_r / \dim \gz_r)}\right )\Pf$.
\end{lemma}

At this point it is convenient to introduce some notation and definitions.
\addtocounter{theorem}{1}
\begin{definition} {\rm The algebra $\gs$ is the {\sl quasi--center} of
$\gn$.  The polynomial function
$\Det_{\gs^*}(\lambda):=\prod_r (\beta_r(\lambda))^{\dim \gg_{\beta_r}}$
on $\gs^*$ is the {\sl quasi--center determinant}.
}
\end{definition}

For $\xi \in \ga$ and $a = \exp(\xi) \in A$ compute 
$
(\Ad(a)\Det_{\gs^*})(\lambda) = \Det_{\gs^*}(\Ad^*(a^{-1})(\lambda))\\
= {\prod}_r (\beta_r(\Ad(a^{-1})^*\lambda))^{\dim \gg_{\beta_r}}
= {\prod}_r (\beta_r(\exp(\beta_r(\xi))\lambda))^{\dim \gg_{\beta_r}}.
$
In other words,

\addtocounter{theorem}{1}
\begin{lemma}\label{trace3}
Let $a = \exp(\xi) \in A$.  Then
$\Ad(a)\Det_{\gs^*} =
\left ( \prod_r \exp(\beta_r(\xi))^{\dim \gz_r}\right )\Det_{\gs^*}$
where $\xi = \log a \in \ga$.
\end{lemma}

Combining Lemmas \ref{trace}, \ref{trace2} and \ref{trace3} we have

\addtocounter{theorem}{1}
\begin{proposition}\label{quasi-invariance}
The product $\Pf\cdot\Det_{\gs^*}$ is an $\Ad(MAN)$--semi--invariant
(and thus $\Ad(AN)$--semi--invariant)
polynomial on $\gs^*$ of degree $\frac{1}{2}(\dim \gn + \dim \gs)$ and of
weight equal to the respective modular functions of $Q$ and $AN$.
\end{proposition}

From $\gn = \gv + \gs$ we have $N = VS$ where
$V = \exp(\gv)$ and $S = \exp(\gs)$.  Now define
\begin{equation}\label{defdp}
D: \text{ Fourier transform of } \Pf\cdot\Det_{\gs^*} \text{, acting on 
the  $S$  variable of $N=VS$.}
\end{equation}

\begin{theorem}\label{dp-min-parab}
The operator $D$ of {\rm (\ref{defdp})} is an invertible
self--adjoint differential operator of degree
$\frac{1}{2}(\dim \gn + \dim \gs)$ on $L^2(MAN)$ with dense
domain $\cC(MAN)$, and it is $\Ad(MAN)$-semi-invariant of
weight equal to the modular function $\delta_{MAN}$\,.
In other words $|D|$ is a Dixmier--Puk\' anszky Operator
on $MAN$ with domain equal to the space of rapidly decreasing $C^\infty$
functions.  This applies as well to $AN$.
\end{theorem}

Since $\lambda \in \gt^*$ has nonzero projection on each summand $\gz_r^*$
of $\gs^*$, and $a \in A$ acts by the positive real scalar
$\exp(\beta_r(\log(a)))$ on $\gz_r$,
\begin{equation}\label{alambda}
A_\lambda = \exp(\{\xi \in \ga \mid \text{each } \beta_r(\xi) = 0\}),
\text{ independent of } \lambda \in \gt^*.
\end{equation}
Because of this independence, and using
$\ga_\diamondsuit = \{\xi \in \ga \mid \text{ each } \beta_r(\xi) = 0\}$,
we define
\begin{equation}\label{adiamond}
A_\diamondsuit = A_\lambda
  \text{ for any (and thus for all) } \lambda \in \gt^*.
\end{equation}

\addtocounter{theorem}{1}
\begin{lemma}\label{ma-diamond}
If $\lambda \in \sigma(\gu^*)$ then the stabilizer
$(MA)_\lambda = M_\diamondsuit A_\diamondsuit$\,.
\end{lemma}

There is no problem with the Mackey obstruction:

\addtocounter{theorem}{1}
\begin{lemma}\label{no-ma-obstruction}
Let $\lambda \in \sigma(\gu^*)$. Recall the extension {\rm (before
(\ref{nm-lambda-family}))}
$\pi_\lambda^\dagger$ of $\pi_\lambda$ to $NM_\diamondsuit$\,.
Then $\pi_\lambda^\dagger$
extends to 
$\widetilde{\pi_\lambda} \in \widehat{NM_\diamondsuit A_\diamondsuit}$
with the same representation space as $\pi_\lambda$\,.
\end{lemma}

When $\lambda \in \sigma(\gu^*)$,
$\widehat{A_\diamondsuit}$ consists of the unitary characters
$\exp(i\phi): a \mapsto e^{i\phi(\log a)}$ with $\phi \in \ga_\diamondsuit^*$.
The representations of $Q$ corresponding to $\lambda$ are the

\begin{equation}\label{lambda-family}
\pi_{\lambda,\gamma,\phi} := \Ind_{NM_\diamondsuit A_\diamondsuit}^{NMA}
  (\widetilde{\pi_\lambda} \otimes \gamma \otimes \exp(i \phi))
  \text{ where } \gamma \in \widehat{M_\diamondsuit} \text{ and }
  \phi \in \ga^*_\diamondsuit\,.
\end{equation} $\Ad^*(A)$ fixes $\gamma$ because $A$ centralizes $M$, and it
fixes $\phi$ because $A$ is commutative, so
\begin{equation} \label{man-equiv}
\pi_{\lambda,\gamma,\phi}\cdot \Ad((ma)^{-1})
 = \pi_{\Ad^*(ma)\lambda, \gamma, \phi}
\end{equation}

\addtocounter{theorem}{1}
\begin{proposition}\label{rep-man}
Plancherel measure for $Q$ is concentrated on the 
the set of all
$\pi_{\lambda,\gamma,\phi}$ for $\lambda \in \sigma(\gu^*)$,
$\gamma \in \widehat{M_\diamondsuit}$\, and
   $\phi \in \ga_\diamondsuit^*$\,.  The equivalence
class of $\pi_{\lambda,\gamma,\phi}$ depends only on
$(\Ad^*(MA)\lambda, \gamma, \phi)$.
\end{proposition}

Representations of $AN$ are the case $\gamma = 1$.  In effect, let
$\pi'_\lambda$ denote the obvious extension $\widetilde{\pi_\lambda}|_{AN}$
of the stepwise square integrable
representation $\pi_\lambda$ from $N$ to $NA_\diamondsuit$
where $\widetilde{\pi_\lambda}$ is given by Lemma \ref{no-ma-obstruction}.
Denote
\begin{equation}\label{irr-na}
        \pi_{\lambda,\phi} = \Ind_{NA_\diamondsuit}^{NA}
        (\pi'_\lambda \otimes \exp(i \phi)) \text{ where }
        \lambda \in \gu^* \text{ and } \phi \in \ga_\diamondsuit^*.
\end{equation} 

\addtocounter{theorem}{1}
\begin{corollary}\label{rep-ma}
Plancherel measure for $AN$ is concentrated on the set of all
$\pi_{\lambda,\phi}$
for $\lambda \in \gu^* \text{ and } \phi \in \ga_\diamondsuit^*$\,.
The equivalence class of $\pi_{\lambda,\phi}$ depends only on
$(\Ad^*(MA)\lambda, \phi)$.
\end{corollary}

A result of C.C. Moore implies 
\addtocounter{theorem}{1}
\begin{lemma}\label{finite-ma-orbit}
The $\Pf$--nonsingular principal orbit set $\gu^*$ is a finite union of
open $\Ad^*(MA)$--orbits.
\end{lemma}

Let $\{\cO_1\,, \dots \cO_v\}$ denote the (open) $\Ad^*(MA)$--orbits on
$\gu^*$.  Denote $\lambda_i = \sigma(\cO_i)$, so 
  $\cO_i = \Ad^*(MA)\lambda_i$ and $(MA)_{\lambda_i} = 
  M_\diamondsuit A_\diamondsuit$ for $1 \leqq i \leqq v.
$
Then Proposition \ref{rep-man} becomes
\begin{theorem}\label{rep-man-ref}
Plancherel measure for $MAN$ is concentrated on
the set $($of equivalence classes of\,$)$
unitary representations
$\pi_{\lambda_i,\gamma,\phi}$ for $1 \leqq i \leqq v$,
$\gamma \in \widehat{M_\diamondsuit}$\, and
   $\phi \in \ga_\diamondsuit^*$\,.
\end{theorem}

Now the Plancherel Theorem for $Q = MAN$ is

The Plancherel Formula (or Fourier Inversion
Formula) for $MAN$ is

\begin{theorem}\label{planch-man}
Let $Q = MAN$ be a minimal parabolic subgroup of the real reductive
Lie group $G$.  Given $\pi_{\lambda,\gamma,\phi} \in \widehat{MAN}$ as
described in {\rm (\ref{lambda-family})} let
$\Theta_{\pi_{\lambda,\gamma,\phi}}: h \mapsto
    \tr \pi_{\lambda,\gamma,\phi}(h)$
denote its distribution character.  Then
$\Theta_{\pi_{\lambda,\gamma,\phi}}$ is a tempered distribution.
If $f \in \cC(MAN)$ then
$$
f(x) = c\sum_{i=1}^{v} \sum_{\gamma \in \widehat{M_\diamondsuit}}
        \int_{\ga^*_\diamondsuit}
        \Theta_{\pi_{\lambda_i,\gamma,\phi}}(D(r(x)f)) |\Pf(b_{\lambda_i})|
        \dim \gamma\,\,d\phi
$$
where $c > 0$ depends on normalizations of Haar measures.
\end{theorem}

The Plancherel Theorem for $NA$ follows similar lines.
For the main computation in the proof of Theorem \ref{planch-man} we
omit $M$ and $\gamma$.  That gives
\begin{equation}\label{one-a-orbit} 
\int_{\ga_\diamondsuit^*}
  \tr \pi_{\lambda_0,\phi}(Dh) \,d\phi
= \int_{\Ad^*(A)\lambda_0} \tr \pi_\lambda(h) |\Pf(b_\lambda)|d\lambda
\end{equation}
In order to go from an $\Ad^*(A)\lambda_0$ to
an integral over $\gu^*$ we use $M$ to parameterize the space of
$\Ad^*(A)$--orbits on $\gu^*$.  If $\lambda \in \gu^*$ one proves
$\Ad^*(A)\lambda \cap \Ad^*(M)\lambda = \{\lambda\}$.
That leads to

\addtocounter{theorem}{1}
\begin{proposition}\label{na-reps}
Plancherel measure for $NA$ is concentrated on the equivalence classes
of representations $\pi_{\lambda,\phi}
= \Ind_{NA_\diamondsuit}^{NA}(\pi'_{\lambda}\otimes \exp(i\phi))$ where
$\lambda \in S_i := \Ad^*(M)\lambda_i$\,, $1 \leqq i \leqq v$,
$\pi'_{\lambda}$ extends
$\pi_\lambda$ from $N$ to $NA_\diamond$ and $\phi \in \ga_\diamond^*$\,.
Representations $\pi_{\lambda,\phi}$ and $\pi_{\lambda',\phi'}$ are
equivalent if and only if $\lambda' \in \Ad^*(A)\lambda$ and
$\phi' = \phi$.  Further, $\pi_{\lambda,\phi}|_N =
\int_{a \in A/A_\diamondsuit} \pi_{\Ad^*(a)\lambda}da$.
\end{proposition}

\begin{theorem}\label{planch-an}
Let $Q = MAN$ be a minimal parabolic subgroup of the real reductive
Lie group $G$.  If $\pi_{\lambda,\phi} \in \widehat{AN}$ let
$\Theta_{\pi_{\lambda, \phi}}: h \mapsto
        \tr \pi_{\lambda,\phi}(h)$
denote its distribution character.  Then
$\Theta_{\pi_{\lambda, \phi}}$ is a tempered distribution.
If $f \in \cC(AN)$ then
$$
f(x) = c\sum_{i=1}^v \int_{\lambda \in \Ad*(M)\lambda_i}
        \int_{\ga^*_\diamondsuit}
        \tr \pi_{\lambda,\phi}(D(r(x)f)) |\Pf(b_\lambda)|
        d\lambda d\phi.
$$
where $c > 0$ depends on normalizations of Haar measures.
\end{theorem}

\section{\hskip -4pt .\hskip 2pt  Parabolic Subgroups in General:
	the Nilradical.}\label{sec9}
\setcounter{lemma}{0}
\setcounter{theorem}{0}
\setcounter{proposition}{0}
\setcounter{corollary}{0}
\setcounter{definition}{0}
\setcounter{remark}{0}
\setcounter{example}{0}
In Sections \ref{sec7} and \ref{sec8} we studied minimal parabolic
subgroups $Q = MAN$ in simple Lie groups, along with certain of their
subgroups $MN$ and $AN$.  This section and the next form a glance at more 
general parabolics.  
This material is taken from \cite{W2015c}, which is a work
in progress, and is limited to the part that I've written down.  
We start with the structure of the nilradical.
\medskip

The condition (c) of (\ref{setup}) does not always hold for nilradicals of
parabolic subgroups.  In this section and the next we weaken (\ref{setup}) to
\begin{equation}\label{setup-weak}
\begin{aligned}
N = &L_1L_2\dots L_{m-1}L_m \text{ where }\\
 &\text{(a) each $L_r$ has unitary representations with coefficients in
$L^2(L_r/Z_r)$,} \\
 &\text{(b) each } N_r := L_1L_2\dots L_r  = N_{r-1}\rtimes L_r 
		\text{ semidirect,}\\
 &\text{(c) if $r \geqq s$ then $[\gl_r,\gz_s] = 0$}.
\end{aligned}
\end{equation}
The conditions of (\ref{setup-weak}) are sufficient
to construct stepwise square integrable representations, but are not always
sufficient to compute the Pfaffian that is the Plancherel density.  So we
refer to (\ref{setup}) as the {\em strong computability condition} and make
make use of the {\em weak computability condition}
\begin{equation}\label{c-weak}
\text{Let }\gl_r = \gl_r' \oplus \gl_r'' \text{ where } 
		\gl_r'' \subset \gz_r \text{ and } \gv_r \subset \gl_r';
\text{ then }[\gl_r,\gl_s] \subset \gl_s'' + \gv_s
        \text{ for } r > s.
\end{equation}
where we retain $\gl_r = \gz_r + \gv_r \text{ and } \gn = \gs + \gv$.
\medskip

Consider an arbitrary parabolic subgroup of $G$.  It contains a 
minimal parabolic $Q = MAN$.  Let $\Psi$ denote
the set of simple roots for the positive system $\Delta^+(\gg,\ga)$.  Then
the parabolic subgroups of $G$ that contain $Q$ are in one to one
correspondence with the subsets $\Phi \subset \Psi$, say $Q_\Phi
\leftrightarrow \Phi$, as follows.  Denote $\Psi = \{\psi_i\}$ and set
\begin{equation}\label{para-roots}
\begin{aligned}
\Phi^{red} &= \left \{\alpha = {\sum}_{\psi_i \in \Psi}
        n_i\psi_i \in \Delta(\gg,\ga) \mid
        n_i = 0 \text{ whenever } \psi_i \notin \Phi \right \} \\
\Phi^{nil} &= \left \{\alpha = {\sum}_{\psi_i \in \Psi}
        n_i\psi_i \in \Delta^+(\gg,\ga) \mid n_i > 0 \text{ for some }
        \psi_i \notin \Phi \right \}.
\end{aligned}
\end{equation}
On the Lie algebra level, $\gq_\Phi = \gm_\Phi + \ga_\Phi + \gn_\Phi$
where
\begin{equation}\label{para-pieces}
\begin{aligned}
&\ga_\Phi = \{ \xi \in \ga \mid \psi(\xi) = 0 \text{ for all } \psi \in
        \Phi \} = \Phi^\perp\,, \\
&\gm_\Phi + \ga_\Phi \text{ is the centralizer of } \ga_\Phi \text{ in }
        \gg \text{, so } \gm_\Phi \text{ has root system } \Phi^{red},
        \text{ and }\\
&\gn_\Phi = {\sum}_{\alpha \in \Phi^{nil}} \gg_\alpha\,, \text{ nilradical
        of } \gq_\Phi\,, \text{ sum of the positive } \ga_\Phi\text{--root
        spaces.}
\end{aligned}
\end{equation}
Since $\gn = \sum_r \gl_r$, as given in (\ref{def-m}) and \ref{def-filtration})
we have
\begin{equation}\label{intersec-1}
\gn_{\Phi} = {\sum}_r (\gn_\Phi \cap \gl_r) =
  {\sum}_r \left ( (\gg_{\beta_r} \cap \gn_\Phi)
        + {\sum}_{\Delta_r^+} (\gg_\alpha \cap \gn_\Phi) \right ).
\end{equation}
As $\ad(\gm)$ is irreducible on each restricted root space,
if $\alpha \in \{\beta_r\} \cup \Delta_r^+$ then
$\gg_\alpha \cap \gn_\Phi$ is $0$ or all of $\gg_\alpha$\,.

\addtocounter{theorem}{1}
\begin{lemma}\label{inter-center}
Suppose $\gg_{\beta_r} \cap \gn_\Phi = 0$.
Then $\gl_r \cap \gn_\Phi = 0$.
\end{lemma}

\addtocounter{theorem}{1}
\begin{lemma}\label{inter-compl}
Suppose $\gg_{\beta_r} \cap \gn_\Phi \ne 0$.  Define $J_r \subset \Delta_r^+$
by $\gl_r \cap \gn_\Phi = \gg_{\beta_r} + \sum_{J_r} \gg_\alpha$\,.  Decompose
$J_r = J'_r \cup J''_r$  where
  $J'_r = \{\alpha \in J_r \mid \sigma_r\alpha \in J_r\}$ and
  $J''_r = \{\alpha \in J_r \mid \sigma_r\alpha \notin J_r\}$.
Then $\gg_{\beta_r} + \sum_{J''_r} \gg_\alpha$ belongs to a single
        $\ga_\Phi$--root space in $\gn_\Phi$\,,  i.e.
        $\alpha|_{\ga_\Phi} = \beta_r|_{\ga_\Phi}$\,, for every
        $\alpha \in J''_r$\,.
\end{lemma}

\addtocounter{theorem}{1}
\begin{lemma} \label{semidirect}
Suppose $\gl_r \cap \gn_\Phi \ne 0$.  Then the algebra
$\gl_r \cap \gn_\Phi$ has center
$\gg_{\beta_r} + \sum_{J''_r} \gg_\alpha$\,, and
$\gl_r \cap \gn_\Phi = (\gg_{\beta_r} + \sum_{J''_r} \gg_\alpha)
+ (\sum_{J'_r} \gg_\alpha))$.
Further, $\gl_r\cap \gn_\Phi
= \left ({\sum}_{J''_r} \gg_\alpha\right ) \oplus
\left ( \gg_{\beta_r} + \left ({\sum}_{J'_r} \gg_\alpha \right )\right )$
direct sum of ideals.
\end{lemma}

It will be convenient to define sets of simple $\ga_\Phi$--roots
\begin{equation}\label{cascade-simple}
\Psi_1 = \Psi \text{ and } \Psi_{s+1} =
        \{\psi \in \Psi \mid  \langle \psi,\beta_i\rangle = 0
        \text{ for } 1 \leqq i \leqq s\}.
\end{equation}
Note that $\Psi_r$ is the simple root system for
$\{\alpha \in \Delta^+(\gg,\ga) \mid \alpha \perp \beta_i
\text{ for } i < r\}$.

\addtocounter{theorem}{1}
\begin{lemma}\label{part-c}
If $r > s$ then
$[\gl_r \cap \gn_\Phi\,,\, \gg_{\beta_s} + {\sum}_{J''_s} \gg_\alpha] = 0$.
\end{lemma}

For our dealings with arbitrary parabolics it is not sufficient to consider 
linear functionals on $\sum_r \gg_{\beta_r}$\,.  Instead we have to look
at linear functionals on
$\sum_r \bigl ( \gg_{\beta_r} + {\sum}_{J''_r} \gg_\alpha \bigr )$.
of the form $\lambda = \sum \lambda_r$ where 
$\lambda_r \in \gg_{\beta_r}^*$ such that $b_{\lambda_r}$ is nondegenerate
on $\sum_r \sum_{J'_r} \gg_\alpha$\,.
We know that (\ref{setup}(c)) holds for the nilradical of the
minimal parabolic $\gq$ that contains $\gq_\Phi$\,. In view of Lemma
\ref{part-c} it follows that
$b_\lambda(\gl_r,\gl_s) = \lambda([\gl_r,\gl_s] = 0$ for $r > s$.
For this particular type of $\lambda$, the bilinear form $b_\lambda$
has kernel $\sum_r \bigl ( \gg_{\beta_s} +
{\sum}_{J''_s} \gg_\alpha \bigr )$ and is nondegenerate on
$\sum_r \sum_{J'_r}\gg_\alpha$\,.  Then
$N_\Phi = (L_1\cap N_\Phi)(L_2\cap N_\Phi)\dots (L_m\cap N_\Phi)$
satisfies the first two conditions of (\ref{setup}).
That is enough to carry out the construction of stepwise square integrable
representations $\pi_\lambda$ of $N_\Phi$\,, but one needs to do more 
to deal with Pfaffian polynomials as in (\ref{setup}(c)) and (\ref{c-weak}).
\medskip

Let $I_1 = \{i \mid \beta_i|_{\ga_\Phi} = \beta_{q_1}|_{\ga_\Phi}\}$ where
$q_1$ is the first index of (\ref{setup}) with
$\beta_{q_1}|_{\ga_\Phi} \ne 0$. 
Next, $I_2 = \{i \mid \beta_i|_{\ga_\Phi} = \beta_{q_2}|_{\ga_\Phi}\}$
where $q_2$ is the first index of (\ref{setup}) such that
$q_2 \notin I_1$ and $\beta_{q_2}|_{\ga_\Phi} \ne 0$.
Continuing as long as possible, 
$I_k = \{i \mid \beta_i|_{\ga_\Phi} = \beta_{q_k}|_{\ga_\Phi}\}$
where $q_k$ is the first index of (\ref{setup}) such that
$q_k \notin (I_1\cup\dots\cup I_{k -1})$ and $\beta_{q_k}|_{\ga_\Phi} \ne 0$.
Then $I_1 \cup \dots \cup I_\ell$ consists of all the
indices $i$ for which $\beta_i|_{\ga_\Phi} \ne 0$.
For $1 \leqq j \leqq \ell$ define
\begin{equation}\label{big-summands}
\gl_{\Phi,j} = {\sum}_{i \in I_j} (\gl_i\cap \gn_\Phi) =
        \Bigl ( {\sum}_{i \in I_j} \gl_i\Bigr ) \cap \gn_\Phi
        \text{ and } \gl^{\dagger}_{\Phi,j} =
        {\sum}_{k \geqq j} \gl_{\Phi,k}\,.
\end{equation}

\addtocounter{theorem}{1}
\begin{lemma}\label{some-brackets}
If $k \geqq j$ then $[\gl_{\Phi,k} , \gl_{\Phi,j}] \subset \gl_{\Phi,j}$\,.
For each index $j$,
$\gl_{\Phi,j}$ and $\gl^{\dagger}_{\Phi,j}$ are subalgebras of $\gn_\Phi$ and
$\gl_{\Phi,j}$ is an ideal in $\gl^{\dagger}_{\Phi,j}$\,.
\end{lemma}

\addtocounter{theorem}{1}
\begin{lemma}\label{not-beta}
If $k > j$ then
$[\gl_{\Phi,k} \,, \gl_{\Phi,j}] \cap \sum_{i \in I_j}\gg_{\beta_i} = 0$.
\end{lemma}

In the notation of Lemma \ref{inter-compl}, if
$r \in I_j$ then
\begin{equation}\label{split-lr}
\gl_r \cap \gn_\Phi = \gl'_r + \gl''_r \text{ where }
        \gl'_r = \gg_{\beta_r} + {\sum}_{J'_r}\gg_\alpha \text{ and }
        \gl''_r = {\sum}_{J''_r}\gg_\alpha\,.
\end{equation}
For $1 \leqq j \leqq \ell$ define
\begin{equation}\label{split-lphi}
\gz_{\Phi,j} = {\sum}_{i \in I_j} (\gg_{\beta_i} +  \gl''_i)
\end{equation}
and decompose
\begin{equation}\label{big-split}
\gl_{\Phi,j} = \gl'_{\Phi,j} + \gl''_{\Phi,j} \text{ where }
        \gl'_{\Phi,j} = {\sum}_{i \in I_j} \gl'_i \text{ and }
        \gl''_{\Phi,j} = {\sum}_{i \in I_j} \gl''_i\,.
\end{equation}

\addtocounter{theorem}{1}
\begin{lemma}\label{central-ideal}
Recall $\gl^{\dagger}_{\Phi,j} = {\sum}_{k \geqq j} \gl_{\Phi,k}$ from
{\rm (\ref{big-summands})}.
For each $j$, both $\gz_{\Phi,j}$ and $\gl''_{\Phi,j}$ are
central ideals in $\gl^{\dagger}_{\Phi,j}$\,, and $\gz_{\Phi,j}$ is the center
of $\gl_{\Phi,j}$.
\end{lemma}

Decompose
\begin{equation}\label{nzv}
\gn_\Phi = \gz_\Phi + \gv_\Phi \text{ where }
\gz_\Phi = \sum_j \gz_{\Phi,j}\,\,,\,\,
\gv_\Phi = \sum_j \gv_{\Phi,j} \text{ and }
\gv_{\Phi,j} = \sum_{i \in I_j} \sum_{\alpha \in J'_i} \gg_\alpha\,.
\end{equation}
Then Lemma \ref{central-ideal} gives us (\ref{c-weak}) for the $\gl_{\Phi,j}$:
$\gl_{\Phi,j} =  \gl'_{\Phi,j} \oplus \gl''_{\Phi,j}$ with
$\gl''_{\Phi,j} \subset \gz_{\Phi,j}$ and $\gv_{\Phi,j} \subset \gl'_{\Phi,j}$.

\addtocounter{theorem}{1}
\begin{lemma} \label{stepwise-nondegen}
For generic $\lambda_j \in \gz_{\Phi,j}^*$ the kernel of $b_{\lambda_j}$ on
$\gl_{\Phi,j}$ is just $\gz_{\Phi,j}$, in other words $b_{\lambda_j}$ is
is nondegenerate on $\gv_{\Phi,j} \simeq \gl_{\Phi,j}/\gz_{\Phi,j}$.
In particular $L_{\Phi,j}$ has square integrable representations.
\end{lemma}

\begin{theorem}\label{gen-setup}
Let $G$ be a real reductive Lie group and $Q$ a real parabolic subgroup.
Express $Q = Q_\Phi$ in the notation of {\rm (\ref{para-roots})} and
{\rm (\ref{para-pieces})}.  Then its nilradical $N_\Phi$ has
decomposition $N_\Phi = L_{\Phi,1}L_{\Phi,2}\dots L_{\Phi,\ell}$
that satisfies the conditions of {\rm (\ref{setup})} and {\rm (\ref{c-weak})}
as follows.
The center $Z_{\Phi,j}$ of $L_{\Phi,j}$ is the analytic subgroup
for $\gz_{\Phi,j}$ and
\begin{equation}\label{gen-setup-list}
\begin{aligned}
 &\text{{\rm (a)} each $L_{\Phi,j}$ has unitary representations
        with coefficients in $L^2(L_{\Phi,j}/Z_{\Phi,j})$,} \\
 &\text{{\rm (b)} each } N_{\Phi,j} :=
        L_{\Phi,1}L_{\Phi,2}\dots L_{\Phi,j}
        \text{ is a normal subgroup of } N_\Phi\\
        &\phantom{XXXXXXXXXXXXXX}
        \text{ with } N_{\Phi,j} = N_{\Phi,j-1}\rtimes L_{\Phi,j}
         \text{ semidirect,}\\
 &\text{{\rm (c)} } [\gl_{\Phi,k},\gz_{\Phi,j}] = 0 \text{ and }
        [\gl_{\Phi,k}, \gl_{\Phi,j}] \subset \gv_{\Phi,j} + \gl_{\Phi,j}''
        \text{ for } k > j.
\end{aligned}
\end{equation}
In particular $N_\Phi$ has stepwise square integrable representations
relative to the decomposition
$N_\Phi = L_{\Phi,1}L_{\Phi,2}\dots L_{\Phi,\ell}$\,.
\end{theorem}

\section{\hskip -4pt .\hskip 2pt  Amenable Subgroups of Semisimple
	Lie Groups.}\label{sec10}
\setcounter{lemma}{0}
\setcounter{theorem}{0}
\setcounter{proposition}{0}
\setcounter{corollary}{0}
\setcounter{definition}{0}
\setcounter{remark}{0}
\setcounter{example}{0}
In this section we apply the results of Section \ref{sec9} to certain
important subgroups of the parabolic $Q_\Phi = M_\Phi A_\Phi N_\Phi$\,,
specifically its amenable subgroups $A_\Phi N_\Phi$\,, $U_\Phi N_\Phi$ and 
$U_\Phi A_\Phi N_\Phi$ where $U_\Phi$ is a maximal compact subgroup of 
$M_\Phi$\,. 
\medskip

The theory of the group $U_\Phi N_\Phi$ goes exactly as in Section \ref{sec7}.
When $N_\Phi = L_{\Phi,1}L_{\Phi,2}\dots L_{\Phi,\ell}$
is weakly invariant we can proceed more or less as in \cite{W2014}.
The argument, but not the final result, will make use of
\addtocounter{theorem}{1}
\begin{definition}\label{invariant}{\rm
The decomposition $N_\Phi =
L_{\Phi,1}L_{\Phi,2}\dots L_{\Phi,\ell}$ of Theorem \ref{gen-setup} is
{\em invariant} if each $\ad(\gm_\Phi)\gz_{\Phi,j} = \gz_{\Phi,j}$,
equivalently if each $\Ad(M_\Phi)\gz_{\Phi,j} = \gz_{\Phi,j}$,
in other words whenever $\gz_{\Phi,j} = \gg_{[\Phi, \beta_{j_0}]}$\,.
The decomposition $N_\Phi = L_{\Phi,1}L_{\Phi,2}\dots L_{\Phi,\ell}$
is {\em weakly invariant} if each $\Ad(U_\Phi)\gz_{\Phi,j} = \gz_{\Phi,j}$\,.
} \hfill $\diamondsuit$
\end{definition}
Set
\begin{equation}\label{regset1}
\gr_\Phi^* = \{\lambda \in \gs_\Phi^* \mid P(\lambda) \ne 0 \text{ and }
        \Ad(U_\Phi)\lambda \text{ is a principal } U_\Phi\text{--orbit on }
        \gs_\Phi^*\}.
\end{equation}
Then $\gr_\Phi^*$ is dense, open and $U_\Phi$--invariant in $\gs_\Phi^*$\,.
By definition of principal orbit the isotropy subgroups of $U_\Phi$ at the
various points of $\gr_\Phi^*$ are conjugate, and we take a measurable
section $\sigma$ to $\gr_\Phi^* \to U_\Phi \backslash \gr_\Phi^*$ on whose
image all the isotropy subgroups are the same,
\begin{equation}\label{m-iso-regset}
U'_\Phi: \text{ isotropy subgroup of } U_\Phi \text{ at }
        \sigma(U_\Phi(\lambda)), \text{ independent of }
        \lambda \in \gr_\Phi^*\,.
\end{equation}
The principal isotropy subgroups
$U'_\Phi$ are pinned down in \cite{HH1970}.
Given $\lambda \in \gr_\Phi^*$ and $\gamma \in \widehat{U'_\Phi}$ let
$\pi_\lambda^\dagger$ denote the
extension of $\pi_\lambda$ to a representation of $U'_\Phi N_\Phi$ on
the space of $\pi_\lambda$ and define 
\begin{equation}\label{rep-un1}
\pi_{\lambda,\gamma} = \Ind_{U'_\Phi N_\Phi}^{U_\Phi N_\Phi} (\gamma \otimes
\pi_\lambda^\dagger).
\end{equation}
The first result in this setting, as in \cite[Proposition 3.3]{W2014}, is

\begin{theorem}\label{planch-mn1}
Suppose that  $N_\Phi = L_{\Phi,1}L_{\Phi,2}\dots L_{\Phi,\ell}$
as in {\rm (\ref{setup-weak})}.  Then the Plancherel
density on $\widehat{U_\Phi N_\Phi}$ is concentrated on the representations
$\pi_{\lambda,\gamma}$ of {\rm (\ref{rep-un1})}, the Plancherel density at
$\pi_{\lambda,\gamma}$ is $(\dim \gamma)|P(\lambda)|$,
and the Plancherel Formula for $U_\Phi N_\Phi$ is
$$
f(un) = c\int_{\gr_\Phi^*/\Ad^*(U_\Phi)}\,
        {\sum}_{\gamma \in \widehat{U'_\Phi}}\,\,
        \tr \Ind_{U'_\Phi N_\Phi}^{U_\Phi N_\Phi}\, r_{un}(f)\cdot
        \dim(\gamma)\cdot |P(\lambda)|d\lambda
$$
where $c = 2^{d_1 + \dots + d_\ell} d_1! d_2! \dots d_\ell!$\, as in
{\rm (\ref{c-d})}.
\end{theorem}

Recall the notion of amenability..  A {\em mean} on a locally compact group 
$H$ is a linear functional $\mu$ on $L^\infty(H)$ of norm $1$ and such that
$\mu(f) \geqq 0$ for all real--valued $f \geqq 0$.  $H$ is {\em amenable}
if it has a left--invariant mean.  
Solvable groups and compact groups are amenable, as are extensions of
amenable groups by amenable subgroups.
In particular $E_\Phi := U_\Phi A_\Phi N_\Phi$ and its closed
subgroups are amenable.
\medskip

We need a technical condition \cite[p. 132]{M1979}.  Let $H$ be the group
of real points in a linear algebraic group whose rational points are Zariski
dense, let $A$ be a maximal $\R$--split torus in $H$, let $Z_H(A)$ denote
the centralizer of $A$ in $H$, and let $H_0$ be the algebraic connected
component of the identity in $H$.  Then $H$ is {\em isotropically connected}
if $H = H_0\cdot Z_H(A)$.  More generally we will say that a subgroup
$H \subset G$ is {\em isotropically connected} if the algebraic hull of
$\Ad_G(H)$ is isotropically connected.  

\begin{proposition}\label{moore}{\rm \cite[Theorem 3.2]{M1979}.}
The groups $E_\Phi := U_\Phi A_\Phi N_\Phi$ are maximal amenable subgroups
of $G$.  They are isotropically connected and self--normalizing.
The various $\Phi \subset \Psi$ are
mutually non--conjugate.  An amenable subgroup $H \subset G$ is contained in
some $E_\Phi$ if and only if it is isotropically connected.
\end{proposition}

The isotropy subgroups are the same at every $\lambda \in \gt_\Phi^*$,
\begin{equation}\label{a-iso-regset}
A'_\Phi: \text{ isotropy subgroup of } A_\Phi \text{ at }
        \lambda \in \gr_\Phi^*\,.
\end{equation}
Given a stepwise square
integrable representation $\pi_\lambda$ where $\lambda \in \gs_\Phi^*$\,,
write $\pi_\lambda^\dagger$ for the extension of $\pi_\lambda$ to a
representation of $A'_\Phi N_\Phi$ on the same Hilbert space.  That
extension exists because the Mackey obstruction vanishes.  The 
representations of $A'_\Phi N_\Phi$
corresponding to $\pi_\lambda$ are the
\begin{equation}\label{rep-an-lambda}
\pi_{\lambda,\phi} := \Ind_{A'_\Phi N_\Phi}^{A_\Phi N_\Phi}
        (\exp(i \phi) \otimes \pi_\lambda^\dagger) \text{ where }
        \phi \in \ga'_\Phi\,.
\end{equation}
Note also that
\begin{equation}\label{rep-an-conj}
\pi_{\lambda,\phi}\cdot \Ad(an) = \pi_{\Ad^*(a)\lambda,\phi}
        \text{ for } a \in A_\Phi \text{ and } n \in N_\Phi\, .
\end{equation}
The resulting formula 
$f(x) = \int_{\widehat{H}} \tr \pi(D(r(x)f))d\mu_H(\pi)$,
$H = A_\Phi N_\Phi$\,, is

\begin{theorem}\label{planch-an1}
Let $Q_\Phi = M_\Phi A_\Phi N_\Phi$ be a parabolic subgroup of the
real reductive Lie group $G$.  Given
$\pi_{\lambda,\phi} \in \widehat{A_\Phi N_\Phi}$ as
described in {\rm (\ref{rep-an-lambda})}, its distribution character
$\Theta_{\pi_{\lambda,\phi}}: h \mapsto
    \tr \pi_{\lambda,\phi}(h)$
is a tempered distribution.
If $f \in \cC(A_\Phi N_\Phi)$ then
$$
f(x) = c \int_{(\ga'_\Phi)^*}
        \left ( \int_{\gs_\Phi^*/\Ad^*(A_\Phi)}
        \Theta_{\pi_{\lambda,\phi}}(D(r(x)f)) |\Pf(b_\lambda)| d\lambda
        \right ) d\phi
$$
where $c = 2^{d_1 + \dots + d_\ell} d_1! d_2! \dots d_\ell!$\,\,.
\end{theorem}

The representations of $U_\Phi A_\Phi N_\Phi$ corresponding to $\pi_\lambda$
are the
\begin{equation}\label{rep-uan-lambda}
\pi_{\lambda,\phi,\gamma} :=\Ind_{U'_\Phi A'_\Phi N_\Phi}^{U_\Phi A_\Phi N_\Phi}
        (\gamma \otimes \exp(i \phi) \otimes \pi_\lambda^\dagger) \text{ where }
        \phi \in \ga'_\Phi \text{ and } \gamma \in \widehat{U'_\Phi}\,.
\end{equation}
Combining Theorems \ref{planch-mn1} and \ref{planch-an1} we arrive at

\begin{theorem}\label{planch-man1}
Let $Q_\Phi = M_\Phi A_\Phi N_\Phi$ be a parabolic subgroup of the
real reductive Lie group $G$ and decompose 
$N_\Phi = L_{\Phi,1}L_{\Phi,2}\dots L_{\Phi,\ell}$
as in {\rm (\ref{setup-weak})}.  Then the Plancherel
density on $\widehat{U_\Phi A_\Phi N_\Phi}$ is concentrated on the 
$\pi_{\lambda,\phi,\gamma}$ of {\rm (\ref{rep-uan-lambda})}, the Plancherel
density at $\pi_{\lambda,\phi,\gamma}$ is $(\dim \gamma)|P(\lambda)|$, the
distribution character $\Theta_{\pi_{\lambda,\phi,\gamma}}: h \mapsto
\tr \pi_{\lambda,\phi,\gamma}(h)$ is tempered, and if
$f \in \cC(U_\Phi A_\Phi N_\Phi)$ then
$$
f(x) = c \sum_{\widehat{U'_\Phi}} \int_{(\ga'_\Phi)^*}
        \left ( \int_{\gs_\Phi^*/\Ad^*(U_\Phi A_\Phi)}
        \Theta_{\pi_{\lambda,\phi,\gamma}}(D(r(x)f)) 
	\deg(\gamma)\,|\Pf(b_\lambda)| d\lambda
        \right ) d\phi
$$
where $c = 2^{d_1 + \dots + d_\ell} d_1! d_2! \dots d_\ell!$\,\,.
\end{theorem}


\begin{thebibliography}{XXXXXXI}

\bibitem{A1963} L. Auslander et al, ``Flows on Homogeneous Spaces'',
Ann. Math. Studies {\bf 53}, 1963.

\bibitem{BB2014} I. Beltita \& D. Beltita, 
Coadjoint orbits of stepwise square integrable representations, to appear.
\{arXiv:1408.1857\}

\bibitem{BL2014} I. Beltita \& J. Ludwig,
Spectral synthesis for coadjoint orbits of nilpotent Lie groups, to appear.
\{arXiv:1412.6323\}

\bibitem{BB2015} I. Beltita \& D. Beltita,
Representations of nilpotent Lie groups via measurable dynamical systems
\{arXiv:1510.05272\}

\bibitem{B1972} 
G. Bredon,
Introduction to Compact Transformation Groups, Academic Press, 1972.

\bibitem{C1989}
W. Casselman,
Introduction to the Schwartz space of $\Gamma \backslash G$,
Canadian J. Math. {\bf 40} (1989), 285--320.

\bibitem{D1972} 
M. Duflo,
Sur les extensions des repr\' esentations irr\' eductibles des groups
de Lie nilpotents,
Ann. Sci. de l' \'Ecole Norm. Sup\' er., 4i\` eme s\' erie {\bf 5} (1972), 
71--120.

\bibitem{DPW2002}
I. Dimitrov, I. Penkov \& J. A. Wolf,
A Bott--Borel--Weil theory for direct limits of algebraic groups,
Amer. J of Math. {\bf 124} (2002), 955--998.

\bibitem{Fa2006}
J. Faraut,
Infinite dimensional harmonic analysis and probability,
in ``Probability Measures on Groups: Recent Directions and Trends,''
ed. S. G. Dani \& P. Graczyk, Narosa, New Delhi, 2006.

\bibitem{G1993}
V. V. Gorbatsevich, A. L. Onishchik \& E. B. Vinberg,
Foundations of Lie Theory and Lie Transformation Groups,
Springer, 1997.

\bibitem{HH1970}
W.-C. Hsiang \& W.-Y. Hsiang,
Differentiable actions of compact connected classical groups II,
Annals of Math. {\bf 92} (1970), 189--223.

\bibitem{K1962} A. A. Kirillov,
Unitary representations of nilpotent Lie groups,
Uspekhi Math. Nauk {\bf 17} (1962), 57--110
(English: Russian Math. Surveys {\bf 17} (1962), 53--104).

\bibitem{Kr1979}
M. Kr\" amer,
Sph\" arische Untergruppen in kompakten zusammenh\" angenden Liegruppen,
Compositio Math. {\bf 38} (1979), 129--153.

\bibitem{LW1978} R. L. Lipsman \& J. A. Wolf,
The Plancherel formula for parabolic subgroups of the classical groups,
Journal D'Analyse Math\' ematique, {\bf 34} (1978), 120--161.

\bibitem{M1979} C. C. Moore,
Amenable subgroups of semi-simple Lie groups and proximal flows.
Israel J. Math. {\bf 34} (1979), 121--138.

\bibitem{M1965} C. C. Moore,
Decomposition of unitary representations defined by discrete subgroups
of nilpotent groups, Ann. Math. {\bf 82} (1965), 146--182

\bibitem{MW1973} C. C. Moore \& J. A. Wolf,
Square integrable representations of nilpotent groups.
Transactions of the American Mathematical Society,
{\bf 185} (1973), 445--462.

\bibitem{N2005}
S. de Neymet Urbina (con la colaboraci\' on de R. Jim\' enez Ben\' itez),
Introducci\' on a los Grupos Topol\' ogicos de Transformaciones,
Sociedad Matem\' atica Mexicana, 2005.

\bibitem{Ol1990}
G. I. Ol'shanskii,
Unitary representations of infinite dimensional pairs $(G,K)$ and the
formalism of R. Howe, in ``Representations of Lie Groups and Related
Topics, ed. A. M. Vershik \& D. P. Zhelobenko,'' Advanced Studies
Contemp. Math. {\bf 7}, Gordon \& Breach, 1990.

\bibitem{P1967} L. Puk\' anszky,
On characters and the Plancherel formula of nilpotent groups,
J. Functional Analysis {\bf 1} (1967), 255--280.

\bibitem{R1972} M. S. Raghunathan,
``Discrete Subgroups of Lie Groups'', Ergebnisse der Mathematik
und ihrer Grenzgebeite {\bf 68}, 1972.

\bibitem{V2001}
E. B. Vinberg,
Commutative homogeneous spaces and co--isotropic symplectic actions,
Russian Math. Surveys {\bf 56} (2001), 1--60.

\bibitem{V2003}
E. B. Vinberg,
Commutative homogeneous spaces of Heisenberg type, Trans Moscow Math. Soc.
{\bf 64} (2003), 45--78.


\bibitem{W1979} J. A. Wolf,
Classification and Fourier inversion for parabolic subgroups with square
integrable nilradical.   Memoirs of the American Mathematical Society,
Number 225, 1979.

\bibitem{W2005}
J. A. Wolf,
Direct limits of principal series representations,  Compositio
Mathematica, {\bf 141} (2005), 1504--1530.

\bibitem{W2007}
J. A. Wolf,
Harmonic Analysis on Commutative Spaces, 
Math. Surveys \& Monographs vol. 142, Amer. Math. Soc., 2007.

\bibitem{W2008}
J. A. Wolf,
Infinite dimensional multiplicity free spaces II:
Limits of commutative nilmanifolds, to appear.

\bibitem{W2013} 
J. A. Wolf,
Stepwise square integrable representations of nilpotent Lie groups,
Mathematische Annalen vol. 357 (2013), pp. 895--914.
\{arXiv: see 1212.1908\}

\bibitem{W2014}
J. A. Wolf, 
The Plancherel Formula for Minimal Parabolic Subgroups, Journal of
Lie Theory, vol. 24 (2014), pp. 791--808.
\{arXiv: 1306.6392 (math RT)\}

\bibitem{W2015a}
J. A. Wolf, 
Stepwise square integrable representations for locally nilpotent
Lie groups, Transformation Groups, vol. 20 (2015), pp. 863--879.
\{arXiv:1402.3828 (math RT, math FA)\}

\bibitem{W2015b}
J. A. Wolf, 
On the analytic structure of commutative nilmanifolds,
The Journal of Geometric Analysis, to appear.
\{arXiv:1407.0399 (math RT, math DG)\}

\bibitem{W2015c}
J. A. Wolf,
Stepwise square integrability for nilradicals of parabolic subgroups 
and maximal amenable subgroups, work in progress.

\bibitem{Y2004}
O. S. Yakimova,
Weakly symmetric riemannian manifolds with reductive isometry group,
Math. USSR Sbornik {\bf 195} (2004), 599--614.

\bibitem{Y2005}
O. S. Yakimova,
``Gelfand Pairs,'' Bonner Math. Schriften (Universit\" at Bonn)
{\bf 374}, 2005.

\bibitem{Y2006}
O. S. Yakimova,
Principal Gelfand pairs,
Transformation Groups {\bf 11} (2006), 305--335.


\end{thebibliography}
\end{document}